\documentclass[preprint]{imsart}

\RequirePackage{amsthm,amsmath,mathrsfs, amssymb, bbm}
\RequirePackage[numbers]{natbib}
\RequirePackage[colorlinks,citecolor=blue,urlcolor=blue]{hyperref}
\RequirePackage{tikz}
\usetikzlibrary{patterns,decorations.pathreplacing}

\newcommand{\prob}{\mathbf{P}}

\newcommand{\exptn}{\mathbf{E}}

\newcommand{\dtv}{\mbox{d}}

\newcommand{\mbbo}{\mathbbm{1}}

\newcommand{\bdelta}{\boldsymbol{\delta}}

\newcommand{\floor}[1]{\lfloor #1 \rfloor}
\newcommand{\iprod}[1]{\langle #1 \rangle}

\newcommand{\beq}{\begin{equation}}
\newcommand{\eeq}{\end{equation}}
\newcommand{\alns}[1]{\begin{align*}#1\end{align*}}
\newcommand{\aln}[1]{\begin{align} #1 \end{align}}
\newcommand{\been}{\begin{enumerate}}
\newcommand{\een}{\end{enumerate}}
\newcommand{\norm}[1]{\| #1 \|}

\newcommand{\inv}{{-1}}

\newcommand{\vconv}{\stackrel{v}{\to}}

\newcommand{\eqd}{\stackrel{d}{=}}

\startlocaldefs
\usepackage{cleveref}
\newtheorem{thm}{Theorem}[section]
\newtheorem{propn}[thm]{Proposition}
\newtheorem{lemma}[thm]{Lemma}
\newtheorem{cor}[thm]{Corollary}

\theoremstyle{remark}
\newtheorem{remark}[thm]{Remark}

\newtheorem{example}{Example}[section]

\theoremstyle{definition}

\endlocaldefs

\allowdisplaybreaks
\numberwithin{equation}{section}

\begin{document}

\begin{frontmatter}
\title{Persistence of heavy-tailed sample averages: principle of infinitely many big jumps}
\runtitle{Persistence of heavy-tailed random walk}

\begin{aug}
\author{\fnms{Ayan} \snm{Bhattacharya} \thanksref{t1, t2}\ead[label=e1]{ayanbhattacharya.isi@gmail.com}},
\author{\fnms{Zbigniew} \snm{Palmowski} \thanksref{t2}\ead[label=e2]{zbigniew.palmowski@pwr.edu.pl}}
\and
\author{\fnms{Bert} \snm{Zwart}\thanksref{t1}\ead[label=e3]{Bert.Zwart@cwi.nl}}

\thankstext{t1}{Partially supported by Dutch Science foundation  NWO VICI grant \#
639.033.413}
\thankstext{t2}{Partially supported by Polish National Science Centre Grant
\#
2018/29/B/ST1/00756 (2019-2022)
}

\runauthor{A. Bhattacharya, Z. Palmowski and B. Zwart}

\affiliation{Centrum Wiskunde \& Informatica\thanksmark{m1} and Wroc\l{a}w University of Science and Technology\thanksmark{m2}}

\address{Centrum Wiskunde \& Informatica\\
P.O. Box 94079\\
1090 GB Amsterdam, Netherlands\\
\printead{e1}
\\
\printead{e3}}

\address{Wroc{\l}aw University of Science and Technology\\
Department of Applied Mathematics\\
Faculty of Pure and Applied Mathematics\\
wyb. Stanis{\l}awa Wyspia\'{n}skiego 27, 50-370 Wroc{\l}aw,
Poland\\
 \printead{e2}\\
 \printead{e3}}
\end{aug}

\begin{abstract}
We consider the sample average of a centered random walk in $\mathbb{R}^d$ with regularly varying step size distribution.
For the first exit time from a compact convex set $A$ not containing the origin, we show that its tail is of lognormal type.
Moreover, we show that the typical way for a large exit time to occur is by having a number of jumps growing logarithmically
in the scaling parameter.
\end{abstract}



\end{frontmatter}

\section{Introduction}

We consider an exit problem for the sample mean of an $\mathbb{R}^d$-valued random walk with zero mean,
where the step size has a distribution which is of multivariate regular variation. Specifically,
let $(\mathbf{X}_i: i\geq 1)$ be an i.i.d.\ sequence of random variables in $\mathbb{R}^d$ ($d\in \mathbb{N}$) such that $\mathbf{X}$
has a multivariate regularly varying distribution with index $\alpha$ (written as $\mathbf{X} \in {\rm RV}(\alpha, \mu)$) where ${\bf X}$ denotes a generic step. Therefore, there exists an increasing  sequence of positive real numbers $(a_n : n \ge 1)$ with $a_n \uparrow \infty$ and a non-null Radon measure $\mu$ on $
\mathscr{B}(\bar{\mathbb{R}}^d \setminus \{\mathbf{0}\})$ with $\mu \big( \bar{\mathbb{R}}^d \setminus \mathbb{R}^d \big) = 0$ such that
\begin{align}
\lim_{n \to \infty} n \prob \big( a_n^\inv \mathbf{X} \in B \big) = \mu (B)   \label{eq:defn:mult:reg:var}
\end{align}
for every  $B \in \mathscr{B}(\bar{\mathbb{R}}^d \setminus \{\mathbf{0}\})$ satisfying $\mu(\partial B ) = 0$ ($\partial B$ denotes the boundary of $B$) and ${\bf 0} \notin \overline{B}$ ($\overline{B}$ denotes the closure of $B$). The limit measure $\mu$ necessarily obeys a homogeneity property, that is, there exists $\alpha \ge 0$ such that $\mu (u \circ B) = u^{-\alpha} \mu(B)$ (where $u \circ B = \{u\cdot \mathbf{x} : \mathbf{x} \in B\}$) for every $u > 0$ and $B \in \mathscr{B} (\bar{\mathbb{R}}^d \setminus \{\mathbf{0}\})$.
We assume that
\begin{equation}\label{alpha}
\alpha >1.
\end{equation}
Additionally, we assume that the $\mathbb{R}^d$-valued random vector ${\bf X}$ satisfies
\begin{equation}\label{zerO}
\exptn \mathbf{X}=\mathbf{0}.
\end{equation}
With $(\mathbf{X}_i: i=1,...,n)$, we associate the random walk
$$\mathbf{S}_k: = \sum_{i=1}^k \mathbf{X}_i,$$
for all $k \in \mathbb{N}.$
In this paper, we investigate the behavior of the survival probability
\begin{equation}
P_{n} := \prob \left(k^\inv \mathbf{S}_k \in A \mbox{ for all } k \in \{1,2, \ldots, n\} \right)
\end{equation}
as $n\rightarrow\infty$, where $A$ is
a compact convex set with non-empty interior that does not contain the origin and
\begin{align}
\prob( X_1 \in A^\circ) >0, \label{eq:ass:positive:prob:jump}
\end{align}  where $A^\circ$ denotes interior of the set $A$. This assumption implies that $P_n >0 $ for every $n$.
On the other hand, (\ref{zerO}) and the LLN subsequently imply
that $P_n\rightarrow 0$ and our aim is to establish its convergence rate. 

Our motivation behind this investigation is two-fold.
First of all, $P_n$ is an example of so-called {\em persistence probability}, that is the probability that
sample average {\it `persists'} in the set $A$ for at least $n$ steps.
It can also be interpreted  as the survival function
$\prob (\tau_A>n)$ of the first time the sample average $\mathbf{S}_k/k$ exits from the set $A$.


Persistence probabilities and related exit problems have recently received a lot of attention in probability theory
and theoretical physics.
In many situations of interest, for a stochastic process in discrete or continuous time and some exit time
$\tau_A$, it turns out that the behavior is either polynomial-like, that is
$\lim_{n\to+\infty} $ $\log \prob(\tau_A>n)/\log n=-\phi$,
or exponential-like, that is
$\lim_{n\to+\infty}\log \prob (\tau_A>n)/n=-\phi$
for a non-negative parameter $\phi$ called the persistence exponent (or survival exponent). This exponent
 usually does not depend on the initial position of the process under consideration.
Random walks and Brownian motions have been analysed in
\cite{MR3342657, doumerc:oconnell:2005, grabiner:1999, 
puchala:rolski:2008, 
wachtel:denisov:2016, vysotsky:2015}.
For results on Gaussian processes, see
\cite{
dembo:mukherjee:2017,
ehrhardt:majumdar:bray:2004, 
wenbo:shao:2004},
and references therein.
If the process under consideration is stationary and one-dimensional, and the set $A$ is a shifted half-line,
 the law of $\tau_A$ corresponds to a first passage time. In this case, fluctuation theory (see \cite{doney:2007}) may be applied; see e.g.\ the survey \cite{aurzada:simon:2015} for an overview concerning mainly
L\'evy processes and (integrated) random walks.
Other one-dimensional processes have been studied; see
for example \cite{hinrichs:kolb:wachtel:2018} for autoregressive sequences.
Recent work on time-homogeneous Markov chains can be found in \cite{aurzada:mukherjee:zeitouni:2017}.
When $\mathbf{E}e^{\langle\mathbf{X}, \lambda\rangle} <\infty$ for all $\lambda \in \mathbb{R}^d$
(hence $\mathbf{X} \notin {\rm RV}(\alpha, \mu)$), the behavior of $P_n$ can be derived
from Mogulskii's theorem, cf.\ \cite[Thm. 5.1.2, p. 176]{DemboZeitouni}.
For a recent survey on persistence probabilities we
refer  to \cite{bray:majumdar:schehr:2013}. 

Our investigation distinguishes from the above-mentioned works by focusing on the sample average $\mathbf{S}_k/k, k\geq 1$,
which is a time-inhomogeneous $\mathbb{R}^d$-valued Markov chain. As mentioned in \cite{bray:majumdar:schehr:2013}, the study of sample averages, and more generally occupation measures, is challenging.
In the case investigated here, we find out that the asymptotics of
$P_n$ is of lognormal type.  That is, there exists a constant $\phi$ depending on the shape of the set $A$ and $\alpha$ such that
\begin{equation}\label{asmain}
\lim_{n\to+\infty} \frac{\log P_n}{(\log n)^2}=-\phi.
\end{equation}
Thus, the behavior of $P_n$ is fundamentally different from the two earlier described cases.
We manage to identify $\phi$ explicitly.
For example, if $d=1$ and $A=[a,b]$  with $0 < a < b$, then the persistence exponent equals
$$\phi=\frac{(\alpha -1)}{2 ( \log b - \log a )}.$$
In the case $d \ge 2$, we provide a simple variational characterization of $\phi$.

An explanation of this untypical asymptotics brings us to our second motivation of this paper,
which is to obtain a sharper understanding of the nature of heavy-tailed large deviations. In turns out that the problem we consider
exhibits a new qualitative phenomenon in the following sense:
we prove that the typical way of getting a large exit time is by having a number of jumps which is growing logarithmic in the scaling parameter $n$.
Hence persistency in our case is caused by
infinitely many large jumps. In other words, the principle of a single big jump used in a significant number of studies
(see \cite{FOSS09anintroduction} and references therein)
does not hold here.

In addition, heavy-tailed sample-path large deviations theorems such as recently
derived in \cite{rhee:blanchet:zwart:2018} do not apply either.
In \cite{rhee:blanchet:zwart:2018},
a sample-path large deviations result for the rescaled random walk $\bar S_n(t), t \in [0,1]$, with $\bar S_n(t) = S_{[nt]}/n$ and $\mathbf{S}_k=S_k$ , has been developed in the case $d=1$.
For a large collection of sets $F$, the results in \cite{rhee:blanchet:zwart:2018} imply that
\begin{equation}
\label{RBZ}
\log \prob \Big( \bar S_n \in F \Big) = - (1+o(1)) J_F (\alpha-1) \log n
\end{equation}
as $n\to+\infty$ with some rate function $J_F$.
This result can be applied to investigate the probability,  for fixed $\epsilon > 0$,
\begin{equation}
P_{\epsilon n,n} := \prob \Big(S_k/k \in [a,b]  \mbox{ for all } k \in \{\lceil \epsilon n\rceil, \ldots, n\} \Big).
\end{equation}
If $- \log \epsilon/ \log (b/a)$ is not an integer, it can be shown that
\begin{equation}\label{eps-asy}
\lim_{n\to+\infty}\frac{\log P_{\epsilon n, n}}{\lceil - \log \epsilon/ \log (b/a) \rceil (\alpha-1) \log n}=1.
\end{equation}
The intuition, which can be made precise using the conditional limit theorems in \cite{rhee:blanchet:zwart:2018}, is that the most likely way for $S_k/k$ to stay in the set $[a,b]$  for $k \in \{\lceil \epsilon n\rceil, \ldots, n\} $ is by having $- \log \epsilon/ \log (b/a)$
large jumps. In the case we are interested in, $O(1)$ jumps will not be sufficient for $S_k/k$ to be persistent.
Therefore, $P_n$ has different asymptotics. Moreover, note that it is tempting to proceed heuristically, and
take $\epsilon = 1/n$ in (\ref{eps-asy}). Apart from not being rigorous, the resulting guess of $\phi$ would actually be off by a factor
$1/2$.

There exist several approaches that can be used to derive the existence, as well as expressions of persistence exponents.
In the case of more general processes,
the Markovian structure is typically exploited. This allows to relate the persistence
exponent to an eigenvalue of an appropriate operator, allowing to marshal analytic methods. This idea is related to
identifying so-called quasi-stationary distributions  (see \cite{MR1465162} for the Brownian motion, \cite{Bogdan, MR3342657, 
MR2248228} for random walks and L\'evy processes,
\cite{MR2986807, MR1334159}
for time-homogeneous Markov processes and \cite{MR3498004, MR2318407, MR2299923} for continuous-time branching processes and Fleming-Viot processes).

\begin{figure}[h]

\begin{tikzpicture}

\draw (0,0) -- (10,5);
\draw (0,0) -- (10,3);
\draw (0,0) -- (10,0);

\draw[thick, red] (2,.6) --(2,1);

\draw(2,0) node{.};
\draw (2,0) node[anchor = north] {$c_1$};

\draw[dashed] (2,1) -- (10/3, 1);

\draw (10/3,0) node{.};
\draw(10/3,0) node[anchor= north]{$\lfloor b c_1/a \rfloor$};

\draw[thick, red] (10/3,1) -- (10/3,5/3);

\draw[dashed] (10/3, 5/3) -- (50/9, 5/3);
\draw (50/9,0) node{.};
\draw(50/9,0) node[anchor = north]{$\lfloor b^2 c_1/a^2 \rfloor $};

\draw[thick, red] (50/9, 5/3) -- (50/9, 25/9);

\draw[dashed] (50/9, 25/9) -- (250/27,25/9);
\draw (250/27,0) node{.};
\draw(250/27,0) node[anchor = north]{$\lfloor b^3 c_1/a^3 \rfloor$};

\draw[thick, red] (250/27, 25/9) -- ( 250/27, 125/27);

\end{tikzpicture}

\caption{ Optimal path for one-dimensional case} \label{figure:optimal:path}

\end{figure}
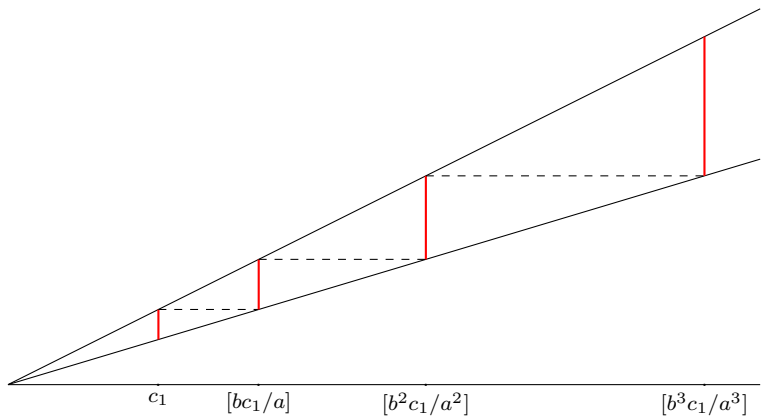

Our work is based on constructing a typical path for the random walk and showing
that this path, sometimes also called the optimal path, is the most likely way for persistence to occur.
For $d=1$ the optimal path is depicted in Figure~\ref{figure:optimal:path} (where the jumps are coloured by red) and it is constructed in the following way.
 Fix a positive finite integer $c_1$. Suppose that the path stays inside the envelope $[ak, bk]$ for all $k \in \{1,2, \ldots , c_1\}$ and the path is at $bc_1$ at time $c_1$. Because of the
 zero drift assumption, the random walk stays around $bc_1$ as long as possible, that is until time $\lfloor bc_1/a \rfloor + 1$. At time $\lfloor bc_1/a \rfloor +1$, it makes the  \textit{first big jump} so that it reaches to the maximum height ($b \lfloor bc_1/a \rfloor + b$) possible and stays there as long as possible, that is until $\lfloor b^2 c_1/ a^2 \rfloor + 1$. Then it again makes a jump. This strategy can be applied recursively, and the resulting path turns out to be the optimal sample path for the event $\{S_k \in [ak, bk] \mbox{ for all } k \in {\mathbb N}\}$. Suppose that $T_i$ denotes the time of the $i$-th jump whose size is denoted by $J_i$. Then we will show that
 a random time $T_i$ can be replaced by $(b/a)^i$ for large enough $i$ with high probability. 
 Let $K_n$ denote the number of big jumps needed until time $n$, i.e. $K_n = \sup \{ i \ge 1: T_i \le n \}$. Then
 $K_n$ can be replaced by  $(\log b/a)^{-1} \log n$ for large $n$ with high probability.
 As we said above, the optimal path can be represented by the random measure
 \begin{equation}
 \sum_{i=1}^{K_n} J_i \bdelta_{T_i},
  \end{equation}
where $\bdelta_x$ is a Dirac measure putting unit mass at $x$.
Moreover,
the probability of a jump of size $J_i$ during $(T_{i-1}, T_i]$ is of order $(b/a)^{i(1- \alpha)}$. Therefore
 $P_n$ is roughly of order $\prod_{i=1}^{\log n} (b/a)^{i(1- \alpha)}$.
 This produces the required estimate
 $\log P_n \asymp -(\alpha -1) (\log b/a)^\inv (\log n)^2/2$ where we write $l(n)\asymp k(n)$ if $\omega_1k(n)\le l(n) \le \omega_2 k(n)$ for some constants $\omega_1$ and $\omega_2$.

\begin{figure}

\begin{tikzpicture}

\draw (0,0) -- (0,5);
\draw (0,0) -- (10,0);

\draw (2,2) -- (5,1);
\draw (5,1) -- (8,3);
\draw (8,3) -- (6,5);
\draw (6,5) -- (3,4);
\draw (3,4) -- (2,2);

\draw[blue] (0,0) -- (8,5);

\draw[dashed] (0,0) -- (8,4);
\draw[dashed] (0,0) -- (6.5,5);


\draw[fill=orange] (16/5,8/5)-- (156/55,24/11) -- (932/143,641/143) --(6.881,3.85) --(16/5,8/5);


\draw (32/5-2,3.5-1.1) node[anchor = south]{$\Gamma(\epsilon)$};

\end{tikzpicture}

\caption{Approximation by constructing a narrow hypercuboid: The blue line denotes the optimal direction. We have constructed the largest possible hypercuboid $\Gamma(\epsilon)$ contained in the intersection of $A$ and the cone. } \label{figure:lower:bound:multidimension}

\end{figure}
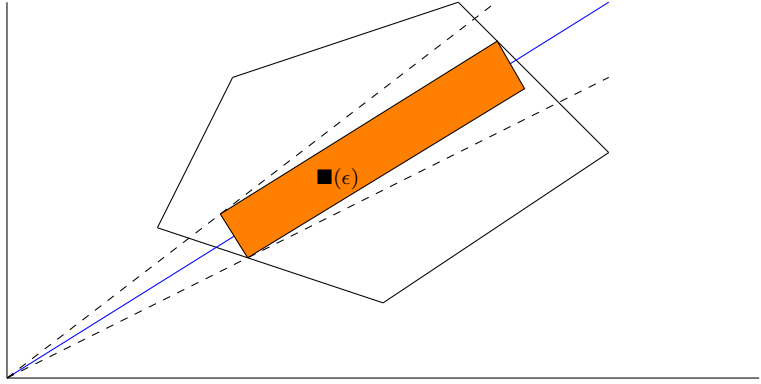

The main idea works also in dimension $d>1$ by choosing an 'optimal' direction $\boldsymbol{\varphi}^*$ that is attaining the supremum
$r^* = \sup_{\boldsymbol{\varphi}\in \Xi(A)} U_{\boldsymbol{\varphi}}/L_{\boldsymbol{\varphi}}$, cf.\ \eqref{r} below.
Using this, we create a convenient inner set of $A$ that is big enough to achieve a sharp enough lower bound for $P_n$.
For this inner set, we take a carefully constructed hypercuboid. A key property is then a certain closure property of a class of hypercuboids
under a direct sum operation.
Another essential feature of our approximation by a sequence of hypercuboids is that we need to allow the fluctuation of the random walk in some directions
though the large jumps happen in the optimal direction $\boldsymbol{\varphi}^*$ only; see Figure \ref{figure:lower:bound:multidimension}.


This paper is organized as follows. In Section \ref{sec:main}, we present the main results Theorem~\ref{thm:persistency:mult:regvar} ($d$-dimensional random walk with $d \ge 2$) and Theorem~\ref{thm:main:theorem:persistence} ($d = 1$) and their consequences along some important examples. In Theorem~\ref{thm:persistency:mult:regvar}, we have assumed the angular measure to be absolutely continuous (with respect to the Lebesgue measure on the surface of the unit sphere) and so, this result does not apply to the $d$-dimensional random walk with independent coordinates (angular measure becomes purely atomic). So in subsection~\ref{subsec:nonstandard}, we present the persistence exponent for a multi-dimensional random walk such that the co-ordinates are independent and the exponent of regular variation might not be the same for every co-ordinate. In section~\ref{sec:proof}, we present the proof of Theorem~\ref{thm:persistency:mult:regvar}. The proof is divided into two parts. In subsection~\ref{subsec:multidimension:upper:bound} and subsection~\ref{subsec:multidimension:lower:bound}, we derive upper and lower bound for the persistence exponent respectively. We further show that the upper and lower bound match and hence, Theorem~\ref{thm:persistency:mult:regvar} follows. The auxiliary results needed to derive the lower bound for the persistence exponent are proved in subsection~\ref{subsec:proof:auxiliary:results}. In subsection~\ref{subsec:proof:persistance:one:dim}, we present a sketch of the proof of Theorem~\ref{thm:main:theorem:persistence}.

\section{Main results}\label{sec:main}

In the definition of regular variation on $\mathbb{R}^d$, we have seen that there exists a Radon measure $\mu$ satisfying the homogeneity property. We first consider $d \ge 2$. The homogeneity property of $\mu$ implies that $\mu$ can also be written as a product measure on $(0, \infty) \times \mathbb{S}^{d-1}$ where $\mathbb{S}^{d-1} = \{\mathbf{x} \in \mathbb{R}^d : \norm{\mathbf{x}} = 1\}$ and $\norm{\mathbf{x}} = \sqrt{x_1^2 + x_2^2 + \ldots + x_d^2 }$. The distance between two sets will be denoted by ${\rm dist}(A, B) = \inf\{\norm{\mathbf{x} - \mathbf{y}} : \mathbf{x} \in A, \mathbf{y} \in B\}$. We need to introduce the polar coordinate transformation to write down the product measure form of $\mu$. The polar co-ordinate transformation is given by $T : \mathbb{R}^{d} \setminus \{\mathbf{0}\} \mapsto (0, \infty) \times \mathbb{S}^{d-1}$, with $T(\mathbf{x}) = (\norm{\mathbf{x}}, \mathbf{x}/\norm{\mathbf{x}})$. This has inverse transformation $T^{\leftarrow} : (0, \infty) \times \mathbb{S}^{d-1} \mapsto \mathbb{R}^d \setminus \{\mathbf{0} \} $ given by $T^{\leftarrow} (r, \mathbf{a}) = r\cdot \mathbf{a}$, where $r \cdot\mathbf{a}$ denotes scalar multiplication of the vector $\mathbf{a}$ and a positive real number $r$. The vector $\mathbf{a}$ can be interpreted as the direction and $r$ is the distance in the direction $\mathbf{a}$.

It is known (e.g.\ Theorem~6.1 in \cite{resnick:2007}) that \eqref{eq:defn:mult:reg:var} is equivalent to  the existence of a Radon measure $\varsigma(\cdot)$ on $\mathbb{S}^{d-1}$ such that
\begin{align}
\lim_{n \to \infty} n \prob \Big( \big(a_n^\inv \norm{\mathbf{X}}, (\norm{\mathbf{X}})^\inv \cdot \mathbf{X} \big) \in C \times D \Big) = \nu_\alpha(C) \varsigma(D),
\end{align}
where $C \in \mathscr{B}((0, \infty))$ and $D \in \mathscr{B}(\mathbb{S}^{d-1})$ and $\nu_\alpha(\cdot)$ is a measure on $(0, \infty)$ such that $\nu_\alpha(x, \infty) = x^{-\alpha}$ for any $x > 0$.  We will assume that the spectral (angular) measure $\varsigma$ is absolutely continuous with respect to the Lebesgue measure on the unit sphere. Note that the spectral measure may not satisfy this assumption: for example it can be atomic if we consider the case where the components of the random vector $\mathbf{X}$ are independent.
Note also that the polar transform is a non-linear transform, that is, the polar transform of a random walk is not a random walk.  Thus, the polar transform can not be used directly to get a one-dimensional positive random walk and compute the persistence exponent from this simpler object. But this decomposition helps to understand the limit. Intuitively, it is clear that the persistence exponent must be based on the radial part of the set under consideration.  

We write $\Xi(B) := \{\norm{\mathbf{x}}^\inv \cdot \mathbf{x} : \mathbf{x} \in B\}$ for any measurable subset $B \in \mathscr{B}(\mathbb{R}^d \setminus \{\mathbf{0}\})$. We consider a 
compact and convex set $A \in \mathscr{B}(\mathbb{R}^d \setminus \{\mathbf{0}\}) $
 which is bounded away from $\mathbf{0}$ ($\mathbf{0} \notin \bar{A}$).
It is clear that $\Xi(A)$ is also compact. We can then write $A = \{r \cdot \boldsymbol{\varphi} :  r \in [L_{\boldsymbol{\varphi}}, U_{\boldsymbol{\varphi}}];~ \boldsymbol{\varphi} \in \Xi(A)\}$ where $L_{\boldsymbol{\varphi}} := \inf \{r : r \cdot \boldsymbol{\varphi} \in A\}$ and $U_{\boldsymbol{\varphi}} : = \sup \{r : r \cdot \boldsymbol{\varphi} \in A\}$. It is clear that $L_{ \boldsymbol{\varphi}}$ and $U_{ \boldsymbol{\varphi}}$ are continuous functions of $\boldsymbol{\varphi}$ as the boundary of a bounded convex set is connected and $L_{\boldsymbol{\varphi}} > 0$ for every $\boldsymbol{\varphi} \in \Xi(A)$ as $A$ is bounded away from $\mathbf{0}$. Thus, we can conclude that $U_{\boldsymbol{\varphi}}/ L_{\boldsymbol{\varphi}}$ is a continuous function of $\boldsymbol{\varphi}$. Define
\begin{align}
r^* := \sup_{\boldsymbol{\varphi} \in \Xi(A)} U_{\boldsymbol{\varphi}}/L_{\boldsymbol{\varphi}} \ge 1. \label{r}
\end{align}
Then there exists ${\boldsymbol{\varphi}}^* \in \Xi(A)$ such that $r^* = U_{\boldsymbol{\varphi}^*}/ L_{\boldsymbol{\varphi}^*}$ as $\Xi(A)$ is compact. This may be non-unique, in which case we fix an arbitrary solution throughout the paper.
Without loss of generality, we can assume that ${\boldsymbol{\varphi}}^*$ points in the direction of the positive orthant of $\mathbb{R}^d$. If it is not the case, then we can rotate the axes to ensure that it holds. We are now ready to present the main result of this work.

\begin{thm} \label{thm:persistency:mult:regvar}
Assume that the angular measure $\varsigma$ is absolutely continuous with respect to the Lebesgue measure, positive on the unit sphere
and the set $A$  with non-empty interior
is compact, convex such that $\mathbf{0} \notin \bar{A}$.
Under the conditions \eqref{alpha}, \eqref{zerO} and \eqref{eq:ass:positive:prob:jump}, we have
\begin{align}
\lim_{n \to \infty} \frac{1}{(\log n)^{2}} \log \prob \Big( k^\inv \mathbf{S}_k \in A \quad\mbox{\rm for all } k =1,2, \ldots, n \Big) = - \frac{\alpha -1}{2(\log r^*)}.
\end{align}
\end{thm}

\begin{remark} 
The persistence exponent $\phi$ and $r^*$ in particular can be computed by developing an alternative representation for $r^*$.
It is not difficult to see that $r^*$ is equal to the largest value of $r$ such that ${\rm dist}(A, r \circ A)=0$ where  $r \circ A = \{r \cdot \mathbf{x}: \mathbf{x} \in A\}$. Since any convex set in $\mathbb{R}^d$ is the intersection of a countable number of half-spaces, there exist vectors $\mathbf{a}_i$ and constants $b_i$ for $i\geq 1$ such that
\begin{equation}
\label{characterization of set A}
A = \{ \mathbf{x}: \langle \mathbf{a}_i, \mathbf{x} \rangle + b_i \leq 0, i\geq 1\}
\end{equation}
where $\langle \mathbf{a}, \mathbf{x} \rangle$ denotes the inner product of vectors $\mathbf{x}$ and $\mathbf{a}$.
Defining the convex function
\begin{equation}
\label{functionH}
H(\mathbf{x}) := \max_i [ \langle \mathbf{a}_i, \mathbf{x}\rangle + b_i],
\end{equation}
the problem of maximizing $r$ such that ${\rm dist} (A,r \circ A)=0$ can now be equivalently written as the solution of the convex program
\begin{equation}
\max_{r,\mathbf{y}} r
\end{equation}
subject to
\begin{equation}
H(\mathbf{y}) \leq 0, H(r \cdot \mathbf{y})\leq 0.
\end{equation}
\end{remark}

\subsection{One-dimensional random walk and interval $[a,b]$}
For $d=1$ and the set $A=[a,b]$ with $0< a < b < \infty$, we consider a collection
$(X_i : i \in \mathbb{N})$ of independent copies of the $\mathbb{R}$-valued, mean-zero
regularly varying random variable $X$ such that
\begin{equation}
\prob \Big( X > x  \Big) = x^{-\alpha} L_+(x) 
 \label{eq:ass:jump:regvar}
\end{equation}
for $x > 0$, such that a tail balance condition
\aln{
\limsup_{x \to \infty} \frac{\prob(X <-x)}{\prob (X > x)} \in [0, \infty)
  \label{eq:ass:jump:heavy:right:tail}
}
holds true, where $L_+$ is a slowly varying function. This is equivalent to assumption \eqref{eq:defn:mult:reg:var} in the case $d=1$.
With $(X_i : i \in \mathbb{N})$, we consider the associated random walk $(S_k : k \ge 1)$ (without using boldface).


\begin{thm}\label{thm:main:theorem:persistence}
Under the assumptions stated above,
\alns{
\lim_{n \to \infty} \frac{1}{( \log n)^2}
\log \prob \Big(  k^\inv S_k \in [a ,b] \quad\mbox{\rm for all } k \in \{1,2, \ldots, n\} \Big) = -\frac{(\alpha -1)}{2 ( \log b - \log a )}
}
for every  $ 0 < a < b < \infty$.
\end{thm}

Note that the above theorem is not a straightforward corollary of Theorem~\ref{thm:persistency:mult:regvar} since the associated angular measure
is necessarily atomic in $d=1$. However, we will briefly show later in the Appendix that its proof follows from the same
steps as the proof of Theorem~\ref{thm:persistency:mult:regvar}.

Theorem~\ref{thm:main:theorem:persistence}
can be used to derive an upper bound for the probability in Theorem~\ref{thm:persistency:mult:regvar} by  projecting a $d$-dimensional random walk in a certain direction. This leads to a natural upper bound
for $P_n$  in terms of a persistence probability for a one-dimensional random walk. In particular, for any $d$-dimensional
vector $\mathbf{c}$,
\begin{equation}
P_n \leq \inf_{{\bf c}: \norm{\bf c} =1} \prob \Big(  k^\inv \langle \mathbf{c}, \mathbf{S}_k \rangle \in \mathbf{c} \bullet A \quad\mbox{\rm for all } k \in \{1,2, \ldots, n\} \Big),
\end{equation}
where
$$ y \in \mathbf{c} \bullet A  \mbox{ if } y=\langle \mathbf{c}, \mathbf{x} \rangle  \mbox{ for some } \mathbf{x} \in A.$$
The assumptions on $A$ and $\mathbf{c}$ imply that $\mathbf{c} \bullet A$ is an interval of the form $[a(\mathbf{c}), b(\mathbf{c})]$.
A natural question is now whether the bound
\begin{equation}
  \phi \geq \sup_{\mathbf{c}: \norm{\mathbf{c}}=1} \frac{(\alpha -1)}{2 ( \log b(\mathbf{c}) - \log a(\mathbf{c}) )}
\end{equation}
for $\phi$ defined in \eqref{asmain} is sharp. This kind of bounding techniques are often applied in light-tailed large deviations.
It can be shown that this bound is sharp if $A$ is a Euclidean ball bounded away from the origin. 
However, if $A$ is a rectangle in the positive
orthant, then the bound is only sharp if and only if the diagonal connecting the southwest corner and northeast corner of
$A$ also passes through the origin. We leave these details as an exercise.

\subsection{Nonstandard regular variation} \label{subsec:nonstandard}
Suppose that  $\mathbf{X} =(X_1, X_2, \ldots, X_d) $ is a random vector such that $X_i$'s are independent and have regularly varying tails with index of regular variation $\alpha_i$ and slowly varying function $L_i(\cdot)$. This is known by the name of nonstandard regular variation in the theory of regular variation (see  \cite[Subsect.~6.5.6]{resnick:2007}).
Then exploiting the independence of components of $\mathbf{S}_k = (S_{k,1}, S_{k,2}, \ldots, S_{k,d}) $ we can get the following easy corollary of Theorem 2.3.
\begin{cor}\label{cor:ind:nonstandard:reg:var}
Suppose that the vector
$\mathbf{X} =(X_1, X_2, \ldots, X_d) $ is such that $X_i$'s are independent and have regularly varying distribution with index of regular variation $\alpha_i$
and each $X_i$ satisfies the assumptions in Theorem~\ref{thm:main:theorem:persistence}.
Then
\begin{align}
& \lim_{n \to \infty} \frac{1}{(\log n)^{2}} \log \prob \left( k^\inv \mathbf{S}_k \in 
\times_{i=1}^d [a_i,b_i],
k =1, \ldots, n \right) \nonumber \\
& \qquad = - \frac{1}{2} \sum_{i=1}^d (\alpha_i -1) (\log b_i - \log a_i)^\inv. \label{nonstandardexponent}
\end{align}
\end{cor}

Note that this cannot be obtained as a corollary of Theorem~\ref{thm:persistency:mult:regvar} as $\mathbf{X} \notin {\rm RV}(\alpha, \mu)$ if the $\alpha_i$'s are not equal. Even if $\alpha_ i = \alpha$ for all $i =1,2, \ldots, d$, then it is known in the literature (see Section~6.5.1 in \cite{resnick:2007}) that the angular measure corresponding to the limit measure $\mu$ is purely atomic and concentrated on the axes which does not fall under the assumptions of Theorem~\ref{thm:persistency:mult:regvar}. Moreover, when all $\alpha$'s are identical, the expression for $\phi$
given in Theorem~\ref{thm:persistency:mult:regvar} does not coincide with the persistence exponent (\ref{nonstandardexponent}). 

\section{Proof of Theorem 2.1} \label{sec:proof}
The proof of  Theorem \ref{thm:persistency:mult:regvar} will be divided into proving the respective asymptotic
lower and upper bounds.

\subsection{Upper bound} \label{subsec:multidimension:upper:bound}

We will show that
\begin{align}
\limsup_{n \to \infty} \frac{1}{(\log n)^2} \log P_n \le - \frac{\alpha - 1}{2 \log r^*}. \label{eq:upper:bound:aim}
\end{align}

\noindent {\bf Step 1.} We divide the set of time points $\{1, 2, \ldots, n\}$ into smaller segments. Fix $\eta > 0$. Then we choose a positive integer $C_1$ such that
\begin{align}
C_1 > 2 + \frac{1}{(1 + \eta) r^* - 1}. \label{eq:choice:C1}
\end{align}
Define $u_0 := 0$, $u_1 := C_1$ and recursively $u_{i+1} := \lfloor (1 +\eta) r^* u_i \rfloor$ for all $i \ge 1$. We also define
\begin{align}
\lambda_n := \sup \{ k > C_1 : u_k \le n\}
\end{align}
for all $n > C_1$. As a consequence, we obtain $u_{\lambda_n} \le n$ and $u_{\lambda_{n} + 1} > n$. Note that $u_{i+1} \ge ( 1+ \eta) r^* u_i - 1$ for all $i \ge 1$. Using these inequalities recursively 
combined with the fact $n \ge u_{\lambda_n}$ yields
\begin{align}
\lambda_n \le 1 +  \frac{\log n}{\log [(1 + \eta)r^*]} - \frac{\log \big[C_1 - (r^* + \eta r^*  - 1)^{-1}\big]}{\log [(1 + \eta)r^*]}. \label{eq:upper:bound:lambdan}
\end{align}
The choice of $C_1$ in \eqref{eq:choice:C1} makes the numerator in the second term in the right hand side of  \eqref{eq:upper:bound:lambdan} well defined. 

 Define $B_i = \{u_{i-1} + 1, u_{i-1} +2, \ldots, u_i\}$ for all $i \ge 1$. Then we have the following bound for $P_n$:
	\begin{align}
	 \prob \Big( \bigcap_{i =1}^n \{ {\bf S}_i \in i \circ A \} \Big) \le \prod_{i =1}^{\lambda_n - 1} \prob \Big( {\bf S}_{u_{i +1}} \in u_{i +1} 		\circ A \Big| \bigcap_{j =1}^i \{{\bf S}_{u_j} \in u_j \circ A\} \Big), \label{eq:upper:bound:prod:cond:prob} 
	\end{align}
using the product formula of conditional probability. 
	
\noindent {\bf Step 2.} Fix $\epsilon_1 \in (0, \eta r^*)$. Then it will be shown in {Step 4} that there exists a positive integer $N(\epsilon_1)$ such that
    \begin{align}
    {\rm dist}(u_i \circ A, u_{i+1} \circ A) \ge u_i C_2 ~~~ \mbox{ for all } i \ge N(\epsilon_1), \label{eq:consecutive:distance}
    \end{align}
where $C_2$ is some positive real number. If $i > N(\epsilon_1)$,  then we can use this property to obtain     \begin{align}
    & \prob \Big( {\bf S}_{u_{i+1}} \in u_{i +1} \circ A;~ {\bf S}_{u_i} \in u_i \circ A; ~\ldots; {\bf S}_{u_1} \in u_1 \circ A \Big) \nonumber \\
    & \le \prob \Big( \norm{ \sum_{j = u_i+ 1}^{u_{i +1}} {\bf X}_j} \ge u_i C_2;~ {\bf S}_{u_i} \in u_i \circ A; \ldots; {\bf S}_{u_1} \in u_1 \circ A \Big) \nonumber \\
    & = \prob \Big( \norm{\sum_{j = u_i +1}^{u_{i +1}} {\bf X}_j } \ge u_i C_2 \Big) \prob \Big( {\bf S}_{u_i} \in u_i \circ A;~ {\bf S}_{u_{i -1}} \in u_{i -1} \circ A;\ldots; {\bf S}_{u_1} \in u_1 \circ A \Big) \label{eq:upper:bound:ind:increment}
    \end{align}
using the independent increment property of the random walk. Combining \eqref{eq:upper:bound:prod:cond:prob} and \eqref{eq:upper:bound:ind:increment}, we infer that
\begin{align}
P_n & \le \prod_{i=1}^{\lambda_n - 1} \prob \big( \norm{{\bf S}_{u_{i+1} - u_i}} > u_i C_2 \big) \nonumber\\
& \le \prod_{i = N(\epsilon_1)+1}^{\lambda_n -1} \prob \Big( (u_{i+1} - u_i)^{-1} \norm{{\bf S}_{u_{i+1} - u_i}} > \frac{1}{(1 + \eta)r^* -1} C_2 \Big),
\end{align}
as $u_i/ (u_{i+1} - u_i)  > [(1+\eta)r^* -1]^{-1}$ for all $i \ge 2$ using $u_{i+1} < u_i (1+\eta)r^*$. Note that $\{(u_{i+1} - u_i)^{-1} \norm{{\bf S}_{u_{i+1} - u_i}} > C_2 [(1+\eta)r^* - 1]^{-1}\} = \{ (u_{i+1} - u_i)^{-1} {\bf S}_{u_{i+1} - u_i} \in \{ {\bf x}: \norm{{\bf x}} > C_2 [(1+\eta)r^* -1]^{-1}\}\}$. This fact leads to the following form of the upper bound for $P_n$:
\begin{align}
\prod_{i = N(\epsilon_1) + 1}^{\lambda_n -1} \prob \Big( (u_{i+1} - u_i)^{-1} {\bf S}_{u_{i+1} - u_i} \in \{{\bf x} : \norm{{\bf x}} > C_2 [(1 + \eta)r^*-1]^{-1}\} \Big).
\end{align}
To bound this expression further, we shall use the following estimate, taken from \cite[Lem.~2.1]{hult:lindskog:mikosch:samorodnitsky:2005}:
\begin{align}
\frac{\prob \Big( n^\inv \mathbf{S}_n \in \cdot \Big)}{n \prob \Big( \norm{\mathbf{X}} > n\Big)} \vconv \mu(\cdot) \label{propn:ldp:mrv:rw:hlms}
\end{align}
on $\mathscr{B}(\mathbb{R}^d \setminus \{\mathbf{0}\})$. Fix $\epsilon_2 > 0$. Note that $(u_{i+1} - u_i) \uparrow \infty$ as $i \to \infty$. So there exists a positive integer $N(\epsilon_2)$ such that
\begin{align}
& \prob \Big( (u_{i+1} - u_i)^{-1} {\bf S}_{u_{i+1} - u_i} \in \{ {\bf x} : \norm{{\bf x}} \ge C_2 [(1 + \eta)r^* -1]^{-1}\} \Big) \nonumber \\
& \le (u_{i+1} - u_i)^{1 - \alpha} L_{\norm{\cdot}} (u_{i+1} - u_i) \big[ \mu \big( \{{\bf x}: \norm{{\bf x}} \ge C_2 [(1 + \eta)r^* -1]^{-1}\} \big) + \epsilon_2 \big] \nonumber \\
& \le C_3 u_i^{1 - \alpha} L_{\norm{\cdot}}(u_i) \label{eq:upper:bound:final:Pn}
\end{align}
for all $i \ge N(\epsilon_2)$, where $L_{\norm{\cdot}}$ is the slowly varying function appearing in the tail distribution function of $\norm{{\bf X}}$, and $C_3$ is some appropriately chosen positive finite real number.  In addition, we have used the fact that $L_{\norm{\cdot}}(u_{i+1} - u_i)/ L_{\norm{\cdot}}(u_i)$ is bounded above as $L_{\norm{\cdot}}$ is a slowly varying function and $u_i^{-1}(u_{i+1} - u_i) \to (1 + \eta)r^* -1 > 0$ as $i \to \infty$.\\

\noindent {\bf Step 3.}  Fix $\epsilon_3 \in (0, \alpha -1)$. We now use  Potter's bound (see e.g. \cite[Prop. 0.8(ii)]{resnick:1987}) which says that there exists an  integer $N(\epsilon_3)$ such that $L_{\norm{\cdot}}(u_i) \le u_i^{\epsilon_3}$ for all $i \ge N(\epsilon_3)$. Define $N_1 = N(\epsilon_1) \vee N(\epsilon_2) \vee N(\epsilon_3)$. Combining the expressions obtained for the upper bound in Step 1 and Step 2, we have
\begin{align}
P_n \le  C_3^{\lambda_n} \prod_{i=1}^{\lambda_n -1} u_i^{1 - \alpha + \epsilon_4}. \label{eq_upper_bound_prodct_formula}
\end{align}
Using the upper bound for $\lambda_n$ obtained in \eqref{eq:upper:bound:lambdan}, straightforward algebra yields
$$P_n \le \exp \Big\{ - \frac{\alpha - 1 - \epsilon_3}{2 \log [(1+\eta)r^*]} (\log n)^2 + O(\log n) \Big\}.$$
The upper bound  \eqref{eq:upper:bound:aim} follows by taking logarithms, dividing by $(\log n)^2$, letting $n\rightarrow\infty$, and finally  $\eta, \epsilon_3 \to 0$. \\

\noindent {\bf Step 4.} Here we shall prove the claim stated in \eqref{eq:consecutive:distance}. We first observe that $u_{i+1}/u_{i} \to (1 + \eta)r^*$ as $i \to \infty$. This implies the existence of a positive integer $N(\epsilon_1)$ such that $u_{i+1}/u_i \ge (1 + \eta)r^* - \epsilon_1$ for all $i \ge N(\epsilon_1)$. We consider $i \ge N(\epsilon_1)$ from now on. It is clear that
$${\rm dist}(u_i \circ A, u_{i+1} \circ A) = u_i \inf_{{\bf x} \in A, {\bf y} \in A} \norm{ (u_i^{-1} u_{i+1}).{\bf x} - {\bf y}} = u_i \inf_{{\bf y} \in A} {\rm dist}({\bf y}, (u_{i+1} u_i^{-1})\circ A).$$
Note that ${\rm dist}({\bf y}, (u_{i+1} u_i^{-1}) \circ A)$ is uniformly continuous in ${\bf y}$. Using that $A$ is compact and every continuous function attains its extrema on a compact set, we conclude that there exists an element ${\bf y}_0 \in A$ such that $\inf_{{\bf y} \in A}{\rm dist}({\bf y}, (u_i^{-1} u_{i+1}) \circ A) = {\rm dist}({\bf y}_0, (u_{i+1}u_i^{-1})\circ A)$. Using continuity of the distance function with the compactness of $A$ once again, we get
$$ \inf_{{\bf y} \in A} {\rm dist}({\bf y}, (u_{i+1} u_i^{-1})\circ A) = \norm{(u_{i+1} u_i^{-1}). {\bf x}_0 - {\bf y}_0}$$
for some pair of elements ${\bf x}_0, {\bf y}_0 \in A$. To prove our claim, it is enough to show that there does not exist any pair of elements ${\bf x}_0$ and ${\bf y}_0$ such that $\norm{(u_{i+1} u_i^{-1}). {\bf x}_0 - {\bf y}_0} = 0$. We prove this by contradiction, so we first assume that there exists a pair ${\bf x}_0$ and ${\bf y}_0$ of elements in $A$ such that $\norm{(u_{i+1} u_i^{-1}). {\bf x}_0 - {\bf y}_0} = 0$ holds. It follows from the property of the Euclidean norm that ${\bf y}_0 = u_{i+1} u_i^{-1}. {\bf x}_0$ and so ${\bf x}_0$ and ${\bf y}_0$ are the vectors in $A$ in the same direction with $\norm{{\bf y}_0}/\norm{{\bf x}_0} = u_{i+1}/u_i \ge r^* + r^* \eta - \epsilon_1 > r^*$. This contradicts the definition of $r^*$ (see \eqref{r}) as ${\bf x}_0, ~{\bf y}_0 \in A$. Hence, the proof is complete.

\subsection{Lower bound} \label{subsec:multidimension:lower:bound}

The proof of the lower bound
\begin{align}
\liminf_{n \to \infty} \frac{1}{(\log n)^{2}} \log P_n \ge - \frac{\alpha - 1}{2\log r^*} \label{eq:upper:bound:eight}
\end{align}
is much more demanding.
Using \eqref{characterization of set A}, and the discussion following that equation, we define $r_\delta$ as the solution of
\begin{equation}
\max_{r,\mathbf{y}} r
\end{equation}
subject to
\begin{equation}
H(\mathbf{y}) \leq \delta, H(r \cdot \mathbf{y})\leq \delta.
\end{equation}
We can equivalently write this as  as the solution of the problem
\begin{equation}
v(\delta) = -r_\delta = \min_{r,\mathbf{y}} -r
\end{equation}
subject to the constraints
\begin{equation}\label{Constr}
H(\mathbf{y}) \leq \delta\quad \text{and}\quad H(r\cdot \mathbf{y})\leq \delta.
\end{equation}
Since $A$ is compact, $H$ has compact level sets for levels $\delta \leq 0$.
Since $H$ is continuous on $A$ and $A$ has a non-empty interior, there exists a $\delta<0$ such that the subset $A^{\delta}:= \{\mathbf{x}:  H(\mathbf{x}) \leq \delta\}$ of $A$ is non-empty,
and so we see that $v(\delta) \leq -1<\infty$ on $\delta$ in a neighborhood of $0$.
Since $H(r \cdot \mathbf{y})$ is a composition of convex functions, it is jointly convex on $[0,\infty) \times [0,\infty)^d$.
Thus, we can apply Theorem 4.2(c) of \cite{BonnansShapiro} with $u = -\delta(1,1)$, $f(x) = - r$, $X= [1, \infty) \times {\mathbb R}^d$ and $G(r,y) = (H(y), H(r.y)): [1,\infty] \times \mathbb{R}^d \to (-\infty,0]^2$
%
to conclude that $v(\delta)$ is continuous
in a neighborhood of $0$.

Consider the set $A^{(\delta)}$ with $\delta < 0$. It is clear that $A^{(\delta)}$ is a proper subset of $A$. Note that $A^{(\delta)}$ is a convex and compact set. So there exist an optimal direction $\varphi^{(\delta)}$ and a straight line (in the direction $\varphi^{(\delta)}$ and the ratio of endpoints $r^{(\delta)}$) such that the straight line is contained in $A^{(\delta)}$. It is immediate that $\varphi^{(\delta)} \in \Xi^\circ$. As a consequence of the Theorem 4.2(c) of \cite{BonnansShapiro}, it follows that $r^{(\delta)} \uparrow r$ as $\delta \uparrow 0$. As $r > 1$, there exists a $\delta_0 < 0$ such that $r^{(\delta)} > 1$ for all $\delta \in (\delta_0, 0)$. Let us fix $\delta \in (\delta_0, 0)$. Note that a segment of straight line (with ratio of the endpoints as $r^{(\delta)}$) in the direction $\varphi^{(\delta)}$ lies in the interior of $A$. Therefore, we can construct a hypercuboid inside the set $A$ (aligned in the direction $\varphi^{(\delta)}$ and the ratio of the endpoints is $r^{(\delta)}$ in the direction $\varphi^{(\delta)}$).  As $r^{(\delta)} > 1$, all the constructions needed to obtain the lower bound in \eqref{eq:final:lower:bound} are possible. Using the same steps, we obtain a lower bound as in \eqref{eq:final:lower:bound} with $r$ replaced by $r^{(\delta)}$. As $\delta$ can be chosen to be arbitrarily close to $0$, we can let $\delta \uparrow 0$ and obtain the desired lower bound.

 Without loss of generality, we can now assume that $\boldsymbol{\varphi}^*\in \Xi^\circ(A)$ for
an interior $\Xi^\circ(A)$ of $\Xi(A)$, that is
$\{\mathbf{y} \in \mathbb{S}^{d-1} : \norm{\mathbf{y} - \boldsymbol{\varphi}^*} < \epsilon\} \subset A$ for some $\epsilon >0 $.
Indeed, if $\boldsymbol{\varphi}^*\in \partial \Xi (A)=\overline{\Xi(A)}\setminus\Xi^\circ(A)$, then
one can consider the set $A^\delta$ instead and then take $\delta \downarrow 0$ at the last step of the proof. \\

\noindent{\bf Strategy of the proof:} In the first step, we shall construct a hypercuboid $\Gamma(\epsilon)$ (the direction of any point in the hypercuboid lies in the $\epsilon$-neighbourhood of $\varphi^*$ in $\Xi^\circ$) which is aligned in the direction $\boldsymbol{\varphi}^*$ and contained in $A$. 
We further show that $\Gamma(\epsilon)$ converges (in the sense of convergence of sets) to the section of the straight line in $A$ in the direction $\boldsymbol{\varphi}^*$ as $\epsilon \downarrow 0$. As we are concerned about the lower bound of $P_n$, we replace $A$ by $\Gamma(\epsilon)$. In the second step, we partition the index set $\{1, 2, \ldots, n\}$ into $(D_i : 1 \le i \le \kappa_n)$ subsets as we did in the upper bound. The main difference here is that the partitions depend on the length $r^{(\epsilon)}$ of $\Gamma(\epsilon)$ in the direction $\varphi^*$, and an auxiliary parameter $\rho$ (choice of $\rho$ depends on $r^{(\epsilon)}$). We also construct $\Upsilon_i \subseteq m_i \circ \Gamma(\epsilon)$ such that $\{ {\bf S}_{m_i} \in \Upsilon_i\}$ for large enough $i$ where $m_i$ is the left end point of $D_i$. We then discuss the strategy for the lower bound in Step~3 and realize the strategy in the rest of the proof.\\

\noindent
\textbf{Steps 1.}(Construction of a hypercuboid inside $A$ approximating the chord of $A$ in the direction $\boldsymbol{\varphi}^*$.)
Define $\mathcal{C}_\epsilon (\boldsymbol{\varphi}^*) := \{\mathbf{y} \in \mathbb{S}^{d-1} : \norm{\mathbf{y} - \boldsymbol{\varphi}^*} \le \epsilon \}$ for some $\epsilon > 0$. From the above assumption $\boldsymbol{\varphi}^* \in \Xi^o(A)$, we can fix $\epsilon > 0$ satisfying $\mathcal{C}_\epsilon(\boldsymbol{\varphi}^*) \subset \Xi(A)$.   This implies that the solid cone $\mathbb{C}_\epsilon := \{r \cdot \mathbf{y} : \mathbf{y} \in \mathcal{C}_\epsilon(\boldsymbol{\varphi}^*), ~ r >0\} $ has non-empty intersection with $A$. We shall say a hypercuboid is aligned in the direction $\boldsymbol{\varphi^*}$ if the hypercuboid is specified by the orthogonal set of unit vectors $(\mathbf{e}_{j} : 1 \le j \le d )$ with $\mathbf{e}_1 = \boldsymbol{\varphi}^*$. We define $\Gamma(\epsilon)$ to be the largest hypercuboid contained in $\mathbb{C}_\epsilon \cap A$. It is clear that as $\epsilon \to 0$, $\mathbb{C}_\epsilon$ converges to the straight line $\{r \cdot  \boldsymbol{\varphi}^* : r > 0\}$. Hence it is clear that $\mathbb{C}_\epsilon \cap A$ converges to $\{r \cdot  \boldsymbol{\varphi^*} : r \in [L_{\boldsymbol{\varphi}^*}, U_{\boldsymbol{\varphi}^*}]\}$. These observations can be used to obtain that $\Gamma(\epsilon)$ converges to $\{r \cdot \boldsymbol{\varphi^*} : r \in [L_{\boldsymbol{\varphi}^*}, U_{\boldsymbol{\varphi}^*}]\}$ using the notion of convergence of sets (as defined in \cite[Def.~4.1]{rockafellar:wets:1998}). Also, note that \eqref{eq:ass:positive:prob:jump} implies that there exists an $\epsilon > 0$ such that $\prob ( {\bf X} \in \{{\bf x} \in A : {\rm dist}({\bf x}, \partial A) >\epsilon \} ) > 0$ for some $\epsilon > 0$. To see this, suppose that $\prob({\bf X} \in A^{-\epsilon}) = 0$ where $A^{-\epsilon} =\{{\bf x} \in A : {\rm dist}({\bf x}, \partial A) >\epsilon \} $ for all $\epsilon>0$. Since $A^\circ = \cup_{\epsilon >0}  A^{-\epsilon}$ and the sets are nested, we can apply the monotone convergence theorem
to conclude that $\prob({\bf X} \in A^\circ) = \lim_{\epsilon \downarrow 0} \prob({\bf X} \in A^{-\epsilon}) =0$ which is a contradiction to \eqref{eq:ass:positive:prob:jump}.  So we can choose $\Gamma(\epsilon)$ for small enough $\epsilon$ such that $\prob({\bf X} \in \Gamma(\epsilon)) >0$.

To specify the $d$-dimensional hypercuboid $\Gamma(\epsilon)$, define
\begin{align}
\Gamma(\epsilon) := \Big\{\mathbf{x} : \langle \mathbf{x}, \mathbf{e}_j \rangle \in [\beta_l^{(j)}(\epsilon), \beta_u^{(j)}(\epsilon)] \mbox{ for all } j = 1,2, \ldots, d\Big\}.
\end{align}
Note that
$\Gamma(\epsilon) \subset A$. Moreover,
we have chosen $(\Gamma(\epsilon) : \epsilon >0 )$ in such a way that
\begin{equation}\label{betas2}
\beta_l^{(j)}(\epsilon) \uparrow 0\quad \mbox{and}\quad \beta_u^{(j)}(\epsilon) \downarrow 0\quad \mbox{as}\quad \epsilon \to 0
\quad \mbox{for all $j =2, \ldots, d$,}
\end{equation}
 and
\begin{equation}\label{betas1}
\beta_l^{(1)}(\epsilon) \downarrow L_{\boldsymbol{\varphi}^*}\quad \mbox{and}\quad \beta^{(1)}_u(\epsilon) \uparrow U_{\boldsymbol{\varphi}^*}\quad \mbox{as}\quad
\epsilon \to 0.
\end{equation}

We have
\begin{align}
&\prob \Big( k^\inv \mathbf{S}_k \in A \mbox{ for all } k =1,2,3, \ldots, n \Big)
\nonumber \\&\qquad\qquad
\ge \prob \Big( k^\inv \mathbf{S}_k \in \Gamma(\epsilon) \mbox{ for all } k =1,2, \ldots, n \Big).  \label{eq:mult:dimension:lower:bound:blacksquare}
\end{align}\\

\noindent\textbf{Steps 2.} (Partitioning the index set and construction of $\Upsilon_i \subseteq m_i \circ \Gamma(\epsilon)$.) As in Step 1 of the upper bound, we divide $\{1,2,\ldots, n\}$ into smaller pieces.
Define
\begin{equation}\label{repsilon}
r^{(\epsilon)} := \beta_u^{(1)}(\epsilon) / \beta_l^{(1)}(\epsilon)
\end{equation}
 and note that $r^{(\epsilon)} \uparrow  r^*$ as $\epsilon \to 0$ by \eqref{betas1}. Fix a constant $\varrho > 0$ small enough so that the following inequality holds:
\aln{
( 1 - \varrho)^2 r^{(\epsilon)} > 1.
}
Note that such a $\varrho$ always exists as $r^{(\epsilon)} > 1$ allows us to choose $\varrho < 1 - (r^{(\epsilon)})^{- 1/2}$. Define
\aln{
m_1 := {\rm C}_1;~~ m_i := \floor{{\rm C}_1 \big[  (1 - \varrho) r^{(\epsilon)} \big]^{i -1}} \mbox{ for all } i \ge 2. \label{eq:defn:mi}
}
For large enough $i$, we define
\aln{
\Upsilon_i & := \Big\{ {\bf x} :  \iprod{ {\bf x}, {\bf e}_1} \in [\beta_l^{(1)}(\epsilon) (1 - \varrho) r^{(\epsilon)} m_i + 10 m_i^{1/\alpha_0 + \delta},~ (1 - \varrho/2) \beta_u^{(1)}(\epsilon) m_i]   \nonumber \\
& \hspace{2cm} \mbox{ and } \iprod{{\bf x}, {\bf e}_j} \in [(2/3) \beta_l^{(j)}(\epsilon) \floor{g m_{i -1}}, ~(2/3) \beta_u^{(j)}(\epsilon) \floor{g m_{i -1}}]  \nonumber \\
& \hspace{3cm}\mbox{ for all } j =2, 3, \ldots, d \Big\}, \label{eq:defn:square_i}
}
where
\aln{
&   g := (1 - \varrho/2) ( 1- \varrho) r^{(\epsilon)} \in ( 1, (1 - \varrho) r^{(\epsilon)}). \label{eq:defn:g}
}
It is easy to check that $\Upsilon_i \subset m_i \circ \Gamma(\epsilon)$ for large enough $i$.
We fix  a large integer ${\rm N}$.  Then we have following lower bound for the right hand side of the expression \eqref{eq:mult:dimension:lower:bound:blacksquare}:
\aln{
 & \prob \Big( k^\inv {\bf S}_k \in \Gamma(\epsilon)  \mbox{ for all } k =1,2, \ldots m_{\rm N} - 1; ~ {\bf S}_{m_{\rm N}} \in \Upsilon_{\rm N} \Big) \prod_{i = \rm N}^{\kappa_n - 1} \prob \Big( k^\inv {\bf S}_k \in \Gamma(\epsilon)  \nonumber \\
 & \hspace{1.5cm}  \mbox{ for all } k =m_i + 1,   m_i + 2, \ldots, m_{i + 1} - 1;~~ {\bf S}_{m_{i +1}} \in \Upsilon_{i + 1}  \Big| {\sf G}_i \Big), \label{eq:lower:bound-product-conditional-probability}
}
where
\aln{
& \kappa_n := \inf \{k : m_k \ge  n\},  \label{defn:kappan}
}
and
\aln{
& {\sf G}_i := \bigg\{k^\inv {\bf S}_k \in \Gamma(\epsilon) \mbox{ for all } k =1,2, \ldots, m_{\rm N} - 1\bigg\} \cap \bigg\{ {\bf S}_{m_{\rm N}} \in \Upsilon_{\rm N}\bigg\} \nonumber \\
& \hspace{2cm} \bigcap_{j = \rm N}^{i -1} \bigg\{ k^\inv {\bf S}_k \in \Gamma(\epsilon) \mbox{ for all } k = m_j + 1, m_j + 2,  \nonumber \\
& \hspace{3cm}  \ldots, m_{j + 1} - 1; ~~ {\bf S}_{m_{j +1}} \in \Upsilon_{j + 1} \bigg\}.
}

\bigskip


\noindent{\bf Step~3.}  We choose the number $g \in (1, (1- \varrho)r^{(\epsilon)})$
in \eqref{eq:defn:g} such that $\Upsilon_{i +1} \subset \floor{g m_i} \circ \Gamma(\epsilon)$. We divide the segment $({\bf S}_k : k \in [m_i + 1, m_{i +1} -1])$ into two parts $({\bf S}_k : k \in [m_i + 1, \floor{g m_i}])$  and $({\bf S}_k : k \in [\floor{g m_i} + 1, m_{i +1}])$. The first part of the segment will be allowed to contribute only to the fluctuation of the random walk where the contribution will be at most of order $m_i^{1/\alpha_0 + \delta} = o(m_i)$. We will use the independent increment property of the random walk and the generalized Kolmogorov's inequality (see \eqref{lemma:gen:kolmogorov}) to show that the probability of this event is close to one for large enough $i$. Observe that the distance between the sets $\Upsilon_i$ and $\Upsilon_{i +1}$ is of order $m_i$, which makes a jump of order $m_i$ necessary. This necessary jump will occur in the second part of the segment  and this part will also contribute to the fluctuation. To analyze this segment we will introduce sets $\Gamma_*$  and $\widetilde{\Gamma}_*$ such that $ m_i \circ \widetilde{\Gamma}_* \subset \Gamma_* $.
Due to the choice of $g$, the jump can occur at any time point in the interval $[\floor{g m_i} + 1, m_{i +1}]$. This strategy, combined with the regular variation of $\norm{{\bf X}}$ and absolute continuity of the law of ${\bf X}/\norm{{\bf X}}$ 
with respect to the uniform angular measure, produces the lower bound
for the probability in \eqref{eq:lower:bound-product-conditional-probability} which is roughly of order $ (m_{i +1} - \floor{g m_i}) \prob( \norm{{\bf X}} > m_i) \approx m_i^{1 - \alpha}$. Then we let $\epsilon \to 0$  and $\varrho \to 0$ to get the desired constants matching the constants in the upper bound.

\bigskip

To realize the strategy, we need the additional sets
	\aln{
	& \widetilde{\Upsilon}_i := \Big\{ {\bf x} : \iprod{{\bf x}, {\bf e}_1} \in [ \beta_l^{(1)}(\epsilon)( 1- \varrho) r^{(\epsilon)} m_i + 9 m_i^{1/\alpha_0 + \delta}, \nonumber \\
	& \hspace{1cm} ~(1 - \varrho/2) \beta_u^{(1)}(\epsilon) m_i + m_i^{1/\alpha_0 + \delta}]  \mbox{ and } \iprod{{\bf x}, {\bf e}_j}  \in [(2/3) \beta_l^{(j)}(\epsilon) \floor{g m_{i -1}} - m_i^{1/\alpha_0 + \delta}, \nonumber \\
	& \hspace{1.5cm} ~~ (2/3) \beta_u^{(j)}(\epsilon) \floor{g m_{i -1}} + m_i^{1/\alpha_0 + \delta}] \mbox{ for all } j =2,3, \ldots, d \Big\}, \label{eq:defn:tilde:square} 
	}
	and
	\aln{
	& \widehat{\Upsilon}_i := \Big\{ {\bf x} : \iprod{ {\bf x}, {\bf e}_1} \in \big[ \beta_l^{(1)}(\epsilon) (1 - \varrho) r^{(\epsilon)} m_i + 8 m_i^{1/\alpha_0 + \delta}, \nonumber \\
	& \hspace{1cm}~ (1 - \varrho/2) \beta_u^{(1)}(\epsilon)m_{i +1} - m_i^{1/\alpha_0 + \delta} \big] \mbox{ and } \iprod{{\bf x}, {\bf e}_j} \in   \big[ (2/3) \beta_l^{(j)}(\epsilon) \floor{g m_{i -1}} - 2 m_i^{1/\alpha_0 + \delta}, \nonumber \\
	& \hspace{ 1.5 cm}~(2/3) \beta_u^{(j)}(\epsilon) \floor{g m_i} - m_i^{1/\alpha_0 + \delta} \big] \mbox{ for all } j =2,3, \ldots, d \Big\}. \label{eq:defn:hat:square}
	}
In the following lemma we introduce their basic properties.
Its proof will be given later in the Appendix.	
\begin{lemma}\label{lemma:subsets:firstsegment:second:segment}
For large enough $i$, we have
\aln{
& \widetilde{\Upsilon}_i \subset \bigcap_{j =m_i + 1}^{\floor{gm_i}} \{j \circ \Gamma(\epsilon) \} \label{eq:inclusion:tilde:square}
}
and
\aln{
&  \widehat{\Upsilon}_i \subset \bigcap_{j = \floor{g m_i} + 1}^{m_{i + 1} - 1} \{ j \circ \Gamma(\epsilon) \}. \label{eq:incllusion:hat:square}
}
\end{lemma}	
Using this lemma, we get the following lower bound for the $i$-th conditional probability in \eqref{eq:lower:bound-product-conditional-probability}:
	\aln{
	& \prob \Big( {\bf S}_k \in \widetilde{\Upsilon}_i \mbox{ for all } k  \in \{m_i + 1, m_i + 2, \ldots, \floor{g m_i}\}; ~ {\bf S}_k  \in \widehat{ \Upsilon}_i \mbox{ for all } k \in \{ \floor{gm_i} + 1,  \nonumber \\
	& \hspace{2cm} \floor{g m_i} + 2, \ldots, m_{i+ 1} - 1 \}  \mbox{ and } ~ {\bf S}_{m_{i +1}} \in \Upsilon_{i +1}  \big| {\sf G}_i \Big) \nonumber \\
	&=  \prob \Big( {\bf S}_k \in \widetilde{\Upsilon}_i \mbox{ for all } k \in {\sf D}_i^{(1)}  \Big| {\sf G}_i \Big) \prob \Big( {\bf S}_k \in \widehat{\Upsilon}_i \mbox{ for all } k \in {\sf D}_i^{(2)}; \nonumber \\
	& \hspace{1.5cm} ~{\bf S}_{m_{i +1}} \in \Upsilon_{i +1} \Big| \Big\{ {\bf S}_k \in \widetilde{\Upsilon}_i  \mbox{ for all } k \in {\sf D}_i^{(1)} \Big\} \cap {\sf G}_i \Big) \nonumber \\
	&  := {\rm T}_i^{(1)} \times {\rm T}_i^{(2)} \label{eq:decomposition:two:terms}
	}
for all large enough $i$, where $${\sf D}_i^{(1)} := \{ m_i + 1, m_i + 2, m_i + 3, \ldots, \floor{g m_i}\} $$ and $${\sf D}_i^{(2)} := \{ \floor{g m_i} + 1, \floor{g m_i} + 2, \ldots, m_{i +1} - 1\}.$$  We shall deal with each of these terms separately.\\

\noindent{\bf Term ${\rm T}_i^{(1)}$}. Note that
	\aln{
	{\rm T}_i^{(1)} & = \big[ \prob({\sf G}_i) \big]^\inv \prob \Big(  \Big\{ {\bf S}_k \in \widetilde{\Upsilon}_i \mbox{ for all } k \in {\sf D}_i^{(1)} \mbox{ and } {\bf S}_{m_i} \in \Upsilon_i \Big\} \cap {\sf G}'_i  \Big),
	}
where
\aln{
{\sf G}'_i & := \bigg\{ k^\inv {\bf S}_k \in \Gamma(\epsilon)  \mbox{ for all } k = 1, 2, \ldots, m_{\rm N} - 1 \mbox{ and } {\bf S}_{m_{\rm N}} \in \Upsilon_{\rm N}\bigg\}   \nonumber \\
& \hspace{.5cm}~  \bigcap_{j ={\rm N}}^{i -2} \bigg\{ k^\inv {\bf S}_k \in \Gamma(\epsilon) \mbox{ for all } k = m_j + 1, m_j + 2, \ldots, m_{m_{j + 1} - 1} \mbox{ and } {\bf S}_{m_j} \in \Upsilon_j \bigg\}  \nonumber \\
& \hspace{1cm} \bigcap \bigg\{ k^\inv {\bf S}_k \in \Gamma(\epsilon) \mbox{ for all } k = m_{i -1} + 1, m_{i -1} + 2, \ldots, m_{i} -1 \bigg\}.
}
Moreover, we have
		\aln{
		& \bigcap_{j =1}^d \Big\{ \max_{k \in {\sf D}_i^{(1)}} \iprod{{\bf S}_k - {\bf S}_{m_i}, {\bf e}_j} \le m_i^{1/\alpha_0 + \delta} ~\mbox{ and }~ \min_{k \in {\sf D}_i^{(1)}} \iprod{{\bf S}_k - {\bf S}_{m_i}, {\bf e}_j} > - m_i^{1/\alpha_0 + \delta}  \Big\}  \nonumber \\
		&  \hspace{.5cm} \bigcap \Big\{{\bf S}_{m_i} \in  \Upsilon_i \Big\} \subset \Big\{ {\bf S}_k \in \widetilde{\Upsilon}_i  \mbox{ for all } k \in {\sf D}_i^{(1)} ~\mbox{ and }~ {\bf S}_{m_i} \in \Upsilon_i \Big\}. \label{eq:onsn:inclusion:term1}
		}
Using this inclusion, we obtain
	\aln{
	{\rm T}_i^{(1)}&  \geq [\prob({\sf G}_i)]^\inv \prob \Big( \Big\{ \max_{k \in {\sf D}_i^{(1)}} \iprod{{\bf S}_k - {\bf S}_{m_i}, {\bf e}_j} \le m_i^{1/\alpha_0 + \delta} ~\mbox{ and }~  \nonumber \\
	& \hspace{1cm} \min_{k \in {\sf D}_i^{(1)}} \iprod{{\bf S}_k - {\bf S}_{m_i}, {\bf e}_j}  > - m_i^{- 1/\alpha_0 + \delta} ~\mbox{ for all } ~ j =1,2, \ldots, d \Big\}  \cap {\sf G}_i\Big). \label{eq:lower:bound:term1:inclusion}
	}
From the independent increment property of the random walk we can conclude that
	\aln{
	 & ({\bf S}_k - {\bf S}_{m_i} : k \in {\sf D}_i^{(1)})  \mbox{ is independent of } ( {\bf S}_j : 1 \le j \le m_i) \mbox{ and  has } \nonumber \\
	  & \mbox{ the same distribution as that of } ({\bf S}_k : 1 \le k \le \floor{g m_i} - m_i). \label{indp:icrement:prop:random:walk:term1}
	}
Thus, the lower bound obtained in \eqref{eq:lower:bound:term1:inclusion} equals
\aln{
& \prob \Big( \Big\{ \max_{k \in {\sf D}_i^{(1)}} \iprod{{\bf S}_k - {\bf S}_{m_i}, {\bf e}_j}  \le m_i^{1/\alpha_0 + \delta} ~\mbox{ and }~  \nonumber \\
& \hspace{2cm} \min_{k \in {\sf D}_i^{(1)}} \iprod{{\bf S}_k - {\bf S}_{m_i}, {\bf e}_j} > - m_i^{1/\alpha_0 + \delta} \mbox{ for all } j =1,2, \ldots, d\Big\} \Big) \nonumber \\
& = \prob \Big( \bigcap_{j =1}^d \Big\{ \max_{k \in {\sf D}_i^{(1)}} \iprod{{\bf S}_{k - m_i}, {\bf e}_j} \le m_i^{1/\alpha_0 + \delta} \mbox{ and }  \nonumber \\
& \hspace{3cm} \min_{k \in {\sf D}_i^{(1)}} \iprod{{\bf S}_{k - m_i}, {\bf e}_j} > - m_i^{1/\alpha_0 + \delta} \Big\} \Big).
}	
We shall now use the positivity of the angular measure to show that the projections of the random walk in each of the directions $({\bf e}_i : 1 \le i \le d)$ are one-dimensional random walks with the same asymptotic tail behaviour. Then, the generalized Kolmogorov inequality (stated in \eqref{lemma:gen:kolmogorov}) is used to obtain the required lower bound. We shall mention the lower bound in the next proposition which will be proved in the Appendix.

\begin{propn} \label{propn:fluctuation:rw}
Fix $\delta \in (0,1)$. Under the assumptions in Theorem~\ref{thm:persistency:mult:regvar}, there exists a large integer ${\rm N}_\delta$ such that  for all $i \ge {\rm N}_\delta$, we have
\aln{
& \prob \Big( \bigcap_{j =1}^d \Big\{ \max_{k \in {\sf D}_i^{(1)}} \iprod{{\bf S}_{k - m_i}, {\bf e}_j} \le m_i^{1/\alpha_0 + \delta} \nonumber\\
& \hspace{2cm} \mbox{ and } \min_{k \in {\sf D}_i^{(1)}} \iprod{{\bf S}_{k - m_i}, {\bf e}_j} > - m_i^{1/\alpha_0 + \delta} \Big\} \Big) > ( 1- \delta).
}
\end{propn}
Thus, from Proposition~\ref{propn:fluctuation:rw} it follows that
\aln{
{\rm T}_i^{(1)} > ( 1- \delta) \label{eq:final:lower:bound:term1}
}
for all large enough $i$'s.

\bigskip

\noindent{\bf Term ${\rm T}_i^{(2)}$}.   Note that
	\aln{
	{\rm T}_i^{(2)} = \frac{\prob \Big( \{ {\bf S}_k \in \widehat{\Upsilon}_i  \mbox{ for all } k \in {\sf D}_i^{(2)};~ S_{m_{i + 1}} \in \Upsilon_{i + 1} \} \cap \{ {\bf S}_k \in \widetilde{\Upsilon}_i \mbox{ for all } k \in {\sf D}_i^{(1)}\} \cap {\sf G}_i \Big)}{\prob \Big(\{ {\bf S}_k \in \widetilde{\Upsilon}_i  \mbox{ for all } k \in {\sf D}_i^{(1)} \} \cap {\sf G}_i \Big)}. \label{eq:ratio:cond:prob:term2}
	}
We now define several sets which will be necessary for the rest of the analysis
\aln{
& \varpi_i := \Big\{ {\bf x} : \iprod{{\bf x}, {\bf e}_1} \in \Big[ - m_i^{ 1/\alpha_0 + \delta}, (1 - \varrho/2) \beta_u^{(1)}(\epsilon) (m_{i + 1 } - m_i) - 2 m_i^{1/\alpha_0 + \delta} \Big] \mbox{ and }\nonumber \\
& \hspace{2cm} \iprod{{\bf x}, {\bf e}_j} \in \Big[- m_i^{1/\alpha_0 + \delta}, (2/3) \beta_u^{(j)}(\epsilon) \big( \floor{g m_i} - \floor{g m_{i -1}}  \big) - 2 m_i^{1/\alpha_0 + \delta} \Big] \nonumber \\
& \hspace{3cm} \mbox{ for all } j =2,3, \ldots, d \Big\}, \\
& \Gamma_* := \Big\{ {\bf x} : \iprod{{\bf x}, {\bf e}_1} \in \Big[ \beta_l^{(1)}(\epsilon) ( 1- \varrho) r^{(\epsilon)}(m_{i +1} - m_i ) + 10 m_{i +1}^{1/\alpha_0 + \delta} - 8 m_i^{1/\alpha_0 + \delta}, \nonumber \\
& \hspace{.7cm} ~~ (1 - \varrho/2) \beta_u^{(1)}(\epsilon) (m_{i +1} - m_i ) - 3 m_i^{1/\alpha_0 + \delta}  \Big]  \mbox{ and } \nonumber \\
& \hspace{1cm} \iprod{{\bf x}, {\bf e}_j} \in \Big[ (2/3) \beta_l^{(j)}(\epsilon) \big( \floor{g m_i}  - \floor{ g m_{i -1}}\big)  - 2 m_i^{1/\alpha_0 + \delta}, \nonumber \\
& \hspace{1cm} (2/3) \beta_u^{(j)} \big( \floor{g m_i} - \floor{g m_{i -1}} \big) - 3 m_i^{1/\alpha_0 + \delta} \Big] \mbox{ for all } j =2,3, \ldots, d \Big\}, \\
&  \widetilde{\varpi}_i := \Big\{ {\bf x} : \iprod{{\bf x}, {\bf e}_1} \in \Big[ \beta_l^{(1)}(\epsilon) (1 - \varrho)r^{(\epsilon)} (m_{i+ 1} - m_i) + 10 m_i^{1/\alpha_0 + \delta} - 9 m_i^{1/\alpha_0 + \delta}, \nonumber \\
 & \hspace{.5cm}  ( 1- \varrho/2) \beta_u^{(1)}(\epsilon)(m_{i +1} - m_i) - m_i^{1/\alpha_0 + \delta} \Big] \mbox{ and } \nonumber \\
 & \hspace{.7cm} \iprod{{\bf x}, {\bf e}_j} \in \Big[ (2/3) \beta_l^{(j)}(\epsilon) \big( \floor{g m_i} - \floor{g m_{i -1}} \big) + m_i^{1/\alpha_0 + \delta}, \nonumber \\
& \hspace{1cm} (2/3) \beta_u^{(j)}(\epsilon) \big( \floor{g m_i} - \floor{g m_{i -1}}\big) - m_i^{1/\alpha_0 + \delta} \Big] \mbox{ for all } j =2,3, \ldots, d \Big\}.
}
Note that the numerator of ${\rm T}_i^{(2)}$ given in \eqref{eq:ratio:cond:prob:term2} has the following lower bound
\aln{
& \prob \Big( \Big\{ {\bf S}_k \in \widehat{\Upsilon}_i \mbox{ for all } k \in {\sf D}_i^{(2)}; ~ {\bf S}_{m_{i +1}} \in \Upsilon_{i +1} \Big\} \cap \Big\{ {\bf S}_k \in \widetilde{\Upsilon}_i \mbox{ for all } k \in {\sf D}_i^{(1)} \Big\} \cap {\sf G}_i \Big) \nonumber \\
&  \ge \prob \Big( \Big\{ {\bf S}_k - {\bf S}_{\floor{gm_i}} \in \varpi_i \mbox{ for all } k \in {\sf D}_i^{(2)};~~ {\bf S}_{m_{i + 1}} - {\bf S}_{\floor{g m_i}} \in \widetilde{\varpi}_i  \Big\} \cap \Big\{  {\bf S}_k \in \widetilde{\Upsilon}_i  \nonumber \\
& \hspace{2.5cm} \mbox{ for all } k \in {\sf D}_i^{(1)} \Big\} \cap {\sf G}_i    \Big) \nonumber \\
& = \prob \Big( {\bf S}_k - {\bf S}_{\floor{g m_i}} \in \varpi_i \mbox{ for all } k \in {\sf D}_i^{(2)};~~ {\bf S}_{m_{i + 1}} - {\bf S}_{\floor{g m_i}} \in \widetilde{\varpi}_i  \Big) \nonumber \\
& \hspace{2.5cm} \prob \Big( \Big\{ {\bf S}_k \in \widetilde{\Upsilon}_i \mbox{ for all } k \in {\sf D}_i^{(1)} \Big\} \cap {\sf G}_i  \Big) \nonumber \\
& = \prob \Big( {\bf S}_{k - \floor{g m_i}} \in \varpi_i  \mbox{ for all } k \in {\sf D}_i^{(2)}; ~~{\bf S}_{m_{i + 1} - \floor{g m_i}} \in \widetilde{\varpi}_i \Big) \nonumber \\
& \hspace{2.5cm} \prob \Big( \Big\{ {\bf S}_k \in \widetilde{\Upsilon}_i \mbox{ for all } k \in {\sf D}_i^{(1)} \Big\} \cap {\sf G}_i \Big). \label{eq:lowerbound:numerator:term2}
}
Finally, we can combine \eqref{eq:ratio:cond:prob:term2} and \eqref{eq:lowerbound:numerator:term2} to have
\aln{
{\rm T}_i^{(2)} \ge & \prob \Big( {\bf S}_{k - \floor{g m_i}} \in \varpi_i \mbox{ for all } k =1,2, \ldots, m_{i + 1} - 1 - \floor{gm_i}; \nonumber \\
& \hspace{2cm}  {\bf S}_{m_{i +1} - \floor{g m_i} } \in \widetilde{\varpi}_i \Big). \label{eq:lower:bound:term2:after:indp:increment}
}

We now decompose the event $\{ {\bf S}_{k - \floor{g m_i}} \in \varpi_i  \mbox{ for all } k = 1,2, \ldots, m_{i +1} - \floor{g m_i} -1 \} \cap \{ {\bf S}_{m_{i + 1} - \floor{ g m_i}  } \in \widetilde{\varpi}_i  \}$ into disjoint events $({\sf E}_t : 1 \le t \le m_{i +1} - \floor{g m_i})$, where
\aln{
 {\sf E}_t & = \Big\{ {\bf X}_t \in \Gamma_* \Big\} \cap \Big\{ \bigcap_{j =1}^d \big\{ \max \big[ \max_{ 1 \le k \le t-1} \iprod{{\bf S}_k, {\bf e}_j}, \nonumber \\
 & \hspace{.5cm} \max_{ t + 1 \le k \le m_{i + 1} - \floor{g m_i}} \iprod{ {\bf S}_k - {\bf X}_t, {\bf e}_j }  \big] \le m_i^{1/\alpha_0 + \delta}   \mbox{ and }  \nonumber \\
 & \hspace{1cm} \min \big[ \min_{1 \le k \le t-1} \iprod{{\bf S}_k, {\bf e}_j}, \min_{ t + 1 \le k \le m_{i + 1} - \floor{gm_i}} \iprod{{\bf S}_k - {\bf X}_t , {\bf e}_j} \big]> - m_i^{1/\alpha_0 + \delta}   \big\} \Big\}. \label{eq:defn:event:Et}
}
It follows from the definition of the events $({\sf E}_t : 1 \le t \le m_{i +1} - \floor{gm_i})$ are exchangeable and hence have the same probabilities. So we have
\aln{
{\rm T}_i^{(2)} & \ge \prob \Big( {\bf S}_{k - \floor{g m_i}} \in \varpi_i \mbox{ for all } k \in {\sf D}_i^{(2)} \mbox{ and } {\bf S}_{m_{i +1}} \in \widetilde{\varpi}_i  \Big) \nonumber \\
&\ge \prob \Big( \bigcup_{t =1}^{m_{i +1} - \floor{g m_i}} {\sf E}_t \Big) = \sum_{t =1}^{m_{i +1} - \floor{g m_i}} \prob({\sf E}_t) = \big( m_{i +1} - \floor{g m_i} \big) \prob({\sf E}_1).  \label{eq:term2:decomposition:disjoint:events}
}
We now estimate $\prob({\sf E}_1)$. Note that
\aln{
\prob({\sf E}_1) & = \prob \Big( \{ {\bf X} \in \Gamma_*\} \cap \Big\{ \bigcap_{j =1}^d \big\{ \max_{2 \le k \le m_{i +1} -\floor{gm_i}} \iprod{ {\bf S}_k - {\bf X}_1, {\bf e}_j} \le m_i^{1/\alpha_0 + \delta} \mbox{ and } \nonumber \\
& \hspace{2cm} \min_{ 2 \le k \le m_{i +1} - \floor{g m_i}} \iprod{{\bf S}_k - {\bf X}_1, {\bf e}_j} > - m_i^{1/\alpha_0 + \delta}  \big\} \Big\} \Big) \nonumber \\
& = \prob \Big( {\bf X} \in \Gamma_* \Big) \prob \Big(  \Big\{ \bigcap_{j =1}^d \big\{ \max_{2 \le k \le m_{i +1} -\floor{gm_i}} \iprod{ {\bf S}_k - {\bf X}_1, {\bf e}_j} \le m_i^{1/\alpha_0 + \delta} \mbox{ and } \nonumber \\
& \hspace{2cm}  \min_{ 2 \le k \le m_{i +1} - \floor{g m_i}} \iprod{{\bf S}_k - {\bf X}_1, {\bf e}_j} > - m_i^{1/\alpha_0 + \delta}  \big\} \Big\} \Big) \nonumber \\
& = \prob({\bf X} \in \Gamma_* )  \prob \Big( \Big\{ \bigcap_{j =1}^d \big\{ \max_{ 1\le k \le m_{i +1} - \floor{gm_i} - 1} \iprod{ {\bf S}_k, {\bf e}_j} \le m_i^{ 1/ \alpha_0 + \delta} \mbox{ and } \nonumber \\
& \hspace{2cm} \min_{ 1 \le k \le m_{i +1} - \floor{g m_i}  - 1} \iprod{ {\bf S}_k, {\bf e}_j } > -  m_i^{1/\alpha_0 + \delta} \big\} \Big\} \Big),  \label{eq:estimate:event:Et}
}
using the independent increment property of the random walk. It can easily be derived from  Proposition~\ref{propn:fluctuation:rw} that
\aln{
&  \prob \Big( \Big\{ \bigcap_{j =1}^d \big\{ \max_{ 1\le k \le m_{i +1} - \floor{gm_i} - 1} \iprod{ {\bf S}_k, {\bf e}_j} \le m_i^{ 1/ \alpha_0 + \delta} \mbox{ and } \nonumber \\
 & \hspace{1.5cm}   \min_{ 1 \le k \le m_{i +1} - \floor{g m_i}  - 1} \iprod{ {\bf S}_k, {\bf e}_j } > -  m_i^{1/\alpha_0 + \delta} \big\} \Big\} \Big) \ge ( 1- \delta)
}
for large enough $i$. We are now left with the estimate of the probability $\prob({\bf X} \in \Gamma_*)$. To do that we shall use the fact that the random variable ${\bf X}$ has a regularly varying tail.  It follows easily from  
$m_i^{1/\alpha_0 + \delta} = o(m_i)$ that  $\{{\bf X} \in \Gamma_* \} \supset \{ m_i^\inv {\bf X}  \in \widetilde{\Gamma}_*\}$ for large enough $i$ where
\aln{
\widetilde{\Gamma}_* & := \Big\{ {\bf x}: \iprod{{\bf x}, {\bf e}_1} \in \Big[ \beta_u^{(1)} (\epsilon) ( 1- 5 \varrho/6) r^{(\epsilon)} \big( (1 - \varrho) r^{(\epsilon)} - 1 \big),  \nonumber \\
& ~~ \beta_u^{(1)}(\epsilon) ( 1- 2 \varrho/3) \big( ( 1- \varrho) r^{(\epsilon) }- 1 \big)  \Big] \mbox{ and } \iprod{{\bf x}, {\bf e}_j} \in  \nonumber \\
&  \hspace{.7cm} \Big[ (1/2) \beta_l^{(j)}(\epsilon) g \big[ 1 -  \big( ( 1- \varrho) r(\epsilon) \big)^\inv \big], ~ (2/3) \beta_u^{(j)}(\epsilon) g \big[ 1 - \big( ( 1- \varrho) r^{(\epsilon)} \big)^\inv \big]  \Big] \nonumber \\
& \hspace{3cm}\mbox{ for all } j =2,3, \ldots, d \Big\}.
}
It is easy to see that $\widetilde{\Gamma}_*$ is bounded away from ${\bf 0}$. Using regular variation of the tail of ${\bf X}$, we get that
\aln{
\lim_{i \to \infty} \frac{ \prob \big( m_i^\inv {\bf X} \in \widetilde{\Gamma}_* \big)}{\prob \big( \norm{{\bf X}} > m_i \big)} = \mu(\widetilde{\Gamma}_*) >0.
}
This means that for large enough $i$, we have
\aln{
\prob \Big( m_i^\inv {\bf X} \in \widetilde{\Gamma}_* \Big)  \ge \prob \big( \norm{{\bf X}} \ge m_i \big) \Big( \mu \big( \widetilde{\Gamma}_* \big) - \varepsilon \Big),
}
where $\varepsilon \in (0, \mu(\widetilde{\Gamma}_*))$ is a fixed small number.  Combining  these facts, we get 
\aln{
\prob({\bf X} \in \Gamma_*) \ge \prob \big( \norm{{\bf X}} \ge m_i \big) \Big( \mu \big( \widetilde{\Gamma}_* \big) - \varepsilon \Big), \label{eq:estimate:lower:bound:jump}
}
for large enough $i$.

Combining \eqref{eq:term2:decomposition:disjoint:events}, \eqref{eq:estimate:event:Et} and \eqref{eq:estimate:lower:bound:jump}, for large enough $i$, we have
\aln{
{\rm T}_i^{(2)} &  \ge (m_{i +1} - \floor{g m_i}) \prob(\norm{{\bf X}} > m_i) (1 - \delta) \nonumber \\
& \sim (1 - \delta) \big((1 -\varrho) r^{(\epsilon)} - g \big)m_i  \prob(\norm{{\bf X}} > m_i) .  \label{eq:final:lower:bound:term2}
}

\bigskip

\noindent{ \bf Steps~4}. It follows from \eqref{eq:decomposition:two:terms}, \eqref{eq:final:lower:bound:term1}, and \eqref{eq:final:lower:bound:term2} that the $i$-th conditional probability in \eqref{eq:lower:bound-product-conditional-probability} can be bounded from below by
\aln{
m_i \prob(\norm{{\bf X}} > m_i) \Big[ (1 - \delta)^2 \big( ( 1- \varrho) r^{(\epsilon)} - g \big) \Big].
}
This estimate, combined with the product formula \eqref{eq:lower:bound-product-conditional-probability}, yields the following lower bound for the probability \eqref{eq:mult:dimension:lower:bound:blacksquare}:
\aln{
& \prob \Big( k^\inv {\bf S}_k \in \Gamma(\epsilon) \mbox{ for all } k = 1, 2, \ldots, m_{\rm N} - 1 \mbox{ and } {\bf S}_{m_{\rm N}} \in \Upsilon_{\rm N} \Big)  \nonumber \\
& \hspace{2cm} \Big[ (1 - \delta)^2 \big( ( 1- \varrho) r^{(\epsilon)} - g \big) \Big]^{\kappa_n}\prod_{i =1 + {\rm N}}^{\kappa_n} \Big[ m_i \prob(\norm{{\bf X}} > m_i)\Big] \nonumber \\
& = {\rm constant} \Big(( 1- \delta^2) [( 1- \varrho)r^{(\epsilon)} - g] \Big)^{\kappa_n} \Big( \prod_{i =1}^{\kappa_n} m_i^{ 1- \alpha} L_{\norm{\cdot}}(m_i) \Big).  \label{eq:product:lower:bound}
}
It follows from the definition of $\kappa_n$ in \eqref{defn:kappan} that
\aln{
\kappa_n = \Big( \log[r^{(\epsilon)}( 1 - \varrho)] \Big)^\inv  \log n + O(1). \label{eq:asymptotics:kappan}
}
We  now use Potter's bound to have $L_{\norm{\cdot}}(m_i) \ge m_i^{ - \eta}$ for large enough $i$ where $\eta$ can be chosen to be arbitrarily small but positive. Now, some straightforward algebra combined with the estimate in \eqref{eq:asymptotics:kappan} leads us to the following
\aln{
\liminf_{n \to \infty} \frac{1}{(\log n)^2} \log P_n & \ge \lim_{n \to \infty} \frac{1}{(\log n)^2} \Big[ \kappa_n \big[ ( 1- \delta)^2 \big( (1 - \varrho )r^{(\epsilon)} - 1 \big) - g \big] - \sum_{i =1}^{\kappa_n} \log m_i \Big] \nonumber \\
& = \frac{1 - \alpha - \eta}{2} \big[ \log (1 - \varrho)r^{(\epsilon)} \big] \lim_{n \to \infty} \frac{\kappa_n^2}{(\log n)^2} \nonumber \\
& = \frac{1 - \alpha - \eta}{ 2} \Big[ \log (1 - \varrho)r^{(\epsilon)} \Big]^{-1}.  \label{eq:final:lower:bound}
}
We can now let $\eta \to 0$, $\varrho \to 0$ and $\epsilon \to 0$ to get the desired constant in the right hand side of \eqref{eq:final:lower:bound}, using the continuity of $r^{(\epsilon)}$ in $\epsilon=0$.

\section{Rest of the proofs}

This section is divided into two subsections.
In subsection~\ref{subsec:proof:auxiliary:results}, we shall first prove the auxiliary results mentioned in subsection~\ref{subsec:multidimension:lower:bound} to derive the lower bound \eqref{eq:final:lower:bound}. In Subsection~\ref{subsec:proof:persistance:one:dim}, we provide a sketch of the proof of Theorem~\ref{thm:main:theorem:persistence}.

\subsection{Proofs of auxiliary results}\label{subsec:proof:auxiliary:results}
\begin{proof}[\bf Proof of Lemma \ref{lemma:subsets:firstsegment:second:segment}]
We first prove \eqref{eq:inclusion:tilde:square}. It is clear from the definition of $m_i$ and $\Gamma(\epsilon)$ that $\{(m_i + 1) \circ \Gamma(\epsilon)\} \cap \{\floor{g m_i} \circ \Gamma(\epsilon) \} = \bigcap_{j = m_i + 1}^{\floor{g m_i} } \{j \circ \Gamma(\epsilon) \}$ and it is enough to show that $\widetilde{\Upsilon}_i \subset \{(m_i + 1) \circ \Gamma(\epsilon)\} \cap \{\floor{g m_i} \circ \Gamma(\epsilon) \}$. To establish this,  we  first consider the direction ${\bf e}_1$.  Note that
\aln{
\big \{\iprod{{\bf x}, {\bf e}_1} : {\bf x} \in \{(m_i + 1) \circ \Gamma(\epsilon)\} \cap \{\floor{g m_i} \circ \Gamma(\epsilon) \} \big\} \supset [g m_i \beta_l^{(1)}(\epsilon), ~ m_i \beta_u^{(1)}(\epsilon)].
}
Comparing the interval with the projection of the set $\widetilde{\Upsilon}_i$ along the direction ${\bf e}_1$, it follows from $(1 - \varrho/2) < 1$ that $g m_i \beta_l^{(1)}(\epsilon) < \beta_l^{(1)}(\epsilon) ( 1- \varrho) r^{(\epsilon)} m_i$ for all $i$ and $(1 - \varrho/2) \beta_u^{(1)}(\epsilon) m_i + m_i^{1/\alpha_0 + \delta} < m_i \beta_u^{(1)}(\epsilon)$ for large enough $i$. We now consider the directions ${\bf e}_j$, where $j =2, 3, \ldots, d$. Fix $j$ and note that
\aln{
\big\{\iprod{{\bf x}, {\bf e}_j} : {\bf x} \in \{ (m_i + 1) \circ \Gamma(\epsilon)  \} \cap \{ \floor{g m_i} \circ \Gamma(\epsilon) \} \big\}  \supset [m_i \beta_l^{(j)}(\epsilon), m_i \beta_u^{(j)}(\epsilon)]
}
as $\beta_l^{(j)}(\epsilon) < 0$  and $\beta_u^{(j)}(\epsilon) > 0$ for all $j =1,2, \ldots, d$.  Comparing this interval with the projection of $\widetilde{\Upsilon}_i $ along the direction ${\bf e}_j$, it follows from $\floor{g m_{i -1}} < m_i$ that $(2/3) \beta_l^{(j)}(\epsilon) \floor{g m_{i -1}} - m_i^{1/\alpha_0 + \delta} > m_i \beta_l^{(j)}(\epsilon)$  for large enough $i$ and $(2/3) \beta_u^{(j)}(\epsilon) m_{i +1} < m_i \beta_u^{(1)}(\epsilon)$  for all $i$.  Hence the inclusion in \eqref{eq:inclusion:tilde:square} follows.

\medskip

We proceed with a proof of \eqref{eq:incllusion:hat:square}. As $ \bigcap_{j = \floor{gm_i} + 1}^{m_{i +1}-1} \{ j \circ \Gamma(\epsilon) \} = \{ (\floor{g m_i} + 1) \circ \Gamma(\epsilon)\} \cap \{(m_{i +1} -1) \circ \Gamma\epsilon) \}$, it will be enough to show that $\widehat{\Upsilon}_i \subset \{ (\floor{g m_i} + 1) \circ \Gamma(\epsilon)\} \cap \{(m_{i +1} -1) \circ \Gamma(\epsilon) \}$.
Consider first the direction ${\bf e}_1$ and note that
\aln{
&&\big\{\iprod{{\bf x}, {\bf e}_1}: {\bf x} \in  \{ (\floor{g m_i} + 1) \circ \Gamma(\epsilon)\} \cap \{(m_{i +1} -1) \circ \Gamma(\epsilon) \} \big \}\nonumber\\&&\qquad \supset [m_{i +1} \beta_l^{(1)}(\epsilon), \floor{g m_i} \beta_u^{(1)}(\epsilon)].
}
Moreover, we have that $m_{i +1} \beta_l^{(1)}(\epsilon) < \beta_l^{(1)}(\epsilon) ( 1- \varrho) r^{(\epsilon)} m_i + 8 m_i^{1/\alpha_0 + \delta}$ for all $i$ and $(1 - \varrho/2) \beta_u^{(1)}(\epsilon) m_{i +1} - m_i^{1/\alpha_0 + \delta} < \floor{g m_i} \beta_u^{(1)}(\epsilon)$ for large enough $i$.
This completes the proof of the inclusion in the ${\bf e}_1$ direction.
Fix now $j \in \{2, 3, \ldots, d\}$.  Then
\aln{
&&\big\{\iprod{{\bf x}, {\bf e}_j}: {\bf x} \in  \{ (\floor{g m_i} + 1) \circ \Gamma(\epsilon)\} \cap \{(m_{i +1} -1) \circ \Gamma(\epsilon) \} \big \}\nonumber\\
&&\qquad \supset [ \floor{g m_i} \beta_l^{(j)}(\epsilon), \floor{g m_i} \beta_u^{(j)}(\epsilon) ]
}
as $\beta_l^{(j)}(\epsilon) < 0$ and $\beta_u^{(j)}(\epsilon) > 0$. Note that $(2/3) \beta_l^{(j)}(\epsilon) \floor{g m_{i -1}} - 2 m_i^{1/\alpha_0 + \delta} > \floor{g m_i} \beta_l^{(j)}(\epsilon)$ for large enough $i$ and $(2/3) \beta_u^{(j)}(\epsilon) \floor{g m_i} - m_i^{1/\alpha_0 + \delta} < \beta_u^{(j)}(\epsilon)$ for all  $i$. This completes the proof of the inclusion stated in \eqref{eq:incllusion:hat:square}.

\end{proof}
\begin{proof}[ \bf Proof of Proposition~\ref{propn:fluctuation:rw}]
Note that
\aln{
& \prob \Big( \bigcap_{j =1}^d \Big\{ \max_{k \in {\sf D}_i^{(1)}} \iprod{{\bf S}_{k - m_i}, {\bf e}_j} \le m_i^{1/\alpha_0 + \delta} \mbox{ and } \min_{k \in {\sf D}_i^{(1)}} \iprod{{\bf S}_{k - m_i}, {\bf e}_j} > - m_i^{1/\alpha_0 + \delta} \Big\} \Big) \nonumber \\
& = 1 - \prob \Big( \bigcup_{j =1}^d \Big\{ \max_{k \in {\sf D}_i^{(1)}} \iprod{{\bf S}_{k - m_i}, {\bf e}_j} \le m_i^{1/\alpha_0 + \delta} \mbox{ and } \min_{k \in {\sf D}_i^{(1)}} \iprod{{\bf S}_{k - m_i}, {\bf e}_j} > - m_i^{1/\alpha_0 + \delta} \Big\}^c  \Big) \nonumber \\
& \ge 1 - \sum_{ j=1}^d \prob \Big( \max_{k \in {\sf D}_i^{(1)}} \iprod{{\bf S}_{k - m_i}, {\bf e}_j} > m_i^{1/\alpha_0 + \delta} \mbox{ or } \min_{k \in {\sf D}_i^{(1)}} \iprod{{\bf S}_{k - m_i}, {\bf e}_j} \le  - m_i^{1/\alpha_0 + \delta}\Big) \nonumber \\
& \ge 1 - \sum_{j =1}^d \Big[ \prob \big(  \max_{k \in {\sf D}_i^{(1)}} \iprod{{\bf S}_{k - m_i}, {\bf e}_j} > m_i^{1/\alpha_0 + \delta} \big) + \prob \big( \min_{k \in {\sf D}_i^{(1)}} \iprod{{\bf S}_{k - m_i}, {\bf e}_j} \le  - m_i^{1/\alpha_0 + \delta} \big) \Big]. \label{eq:union:bound:fluctuation:event}
}
Fix $\delta_{j,t}$ for $j=1,2, \ldots d$ and $t =1,2$ such that $\delta_{j,t}  \in (0, \delta/2d)$.
We claim
\begin{equation}\label{maxin}
\prob \Big( \max_{ 1 \le k \le \lfloor g m_i \rfloor - m_i} \langle \mathbf{S}_k, \mathbf{e}_j \rangle  > m_i^{1/\alpha_0 + \delta} \Big) < \delta_{j,1}
\end{equation}
for sufficiently large $i$.
To prove it we will use the following lemma.
\begin{lemma} \label{lemma:directed:rv:random:walk}
Let $\mathbf{X} \in {\rm RV}(\alpha, \mu)$and $\mu = \nu_\alpha \otimes \varsigma$ on $(0, \infty) \times \mathbb{S}^{d-1}$ with $\varsigma$ being absolutely continuous with respect to the Lebesgue measure. Then for any direction vector $\mathbf{u} \in \mathbb{S}^{d-1}$, we have $\langle \mathbf{u}, \mathbf{X} \rangle \in {\rm RV}( \alpha, \vartheta_\alpha)$ where $\vartheta_\alpha$ is a Radon measure on $\mathbb{R} \setminus \{0\}$ with
\begin{align}
& \vartheta_\alpha (dx) := \alpha \mu(\{\mathbf{y}: \langle \mathbf{u}, \mathbf{y} \rangle > 1\}) x^{-\alpha -1} dx  \mbbo(x > 0) \nonumber \\
& \hspace{2cm} + \alpha \mu( \{ \mathbf{y} : \langle \mathbf{u}, \mathbf{y} \rangle < -1 \}) (-x)^{-\alpha -1} \mbbo(x < 0).
\end{align}
\end{lemma}
Using Lemma~\ref{lemma:directed:rv:random:walk}, note that $\langle \mathbf{S}_k, \mathbf{e}_j \rangle=\sum_{i'=1}^k Y_{i'}^{(j)}$ is a mean $0$ random walk with steps
$Y_{i'}^{(j)}=\langle \mathbf{X}_{i'}, \mathbf{e}_j \rangle \in { \rm RV}(\alpha, \vartheta_\alpha)$ for all $j=1,2, \ldots, d$.
For $\alpha \in (1,2]$, we will apply the generalized Kolmogorov inequality given in \cite{shneer:wachtel:2009}:
\begin{align}
\prob \left( \max_{1 \le k \le m} \sum_{i'=1}^k Y_{i'} \ge x \right) \le C_4 m x^{-2} \exptn \left[ (Y^{(j)})^2 \mbbo\left(|Y^{(j)}| < x\right)\right],\label{lemma:gen:kolmogorov}
\end{align}
where $C_4$ is some constant and $Y^{(j)} = \langle {\bf X}, {\bf e}_j \rangle$. In this case, as \cite{shneer:wachtel:2009} noted, $\exptn \left[ (Y^{(j)})^2 \mbbo\left(|Y^{(j)}| < x\right)\right]$ is regularly
varying with index $2-\alpha$ (or slowly varying if $\alpha =2$).
For $\alpha >2$ we can apply the classical Kolmogorov inequality. In both cases we can bound
$$\prob \Big( \max_{ 1 \le k \le \lfloor g m_i \rfloor - m_i} \langle \mathbf{S}_k, \mathbf{e}_j \rangle  > m_i^{1/\alpha_0 + \delta} \Big)\leq
C_5m_i^{-\alpha_0 \delta +\eta},
$$
where $\eta$ appears due to Potter's bound applied to the slowly varying part of $\exptn \left[ (Y^{(j)})^2 \mbbo\left(|Y^{(j)}| < x\right)\right]$
and $C_5$ is some constant.
For $\eta>0$ sufficiently small, this upper bound gives \eqref{maxin} as $m_i\to \infty$ with $i\to\infty$.

Similarly, we can prove that
$$\prob \left(  \min_{ 1 \le k \le \lfloor g m_i \rfloor - m_i} \langle \mathbf{S}_k, \mathbf{e}_j \rangle < - m_i^{1/\alpha_0 + \delta} \right) < \delta_{j,2}$$
for large enough $i$. Hence the proof of the proposition follows from the lower bound obtained in \eqref{eq:union:bound:fluctuation:event}.
\end{proof}

\begin{proof}[\bf Proof of Lemma~\ref{lemma:directed:rv:random:walk}]
To prove this lemma, we need to find $(b_n : n \ge 1)$ such that
\begin{align}
\lim_{n \to \infty} n \prob \Big( b_n^\inv \langle \mathbf{X}, \mathbf{u} \rangle \in B \Big) = \vartheta_\alpha(B) \in (0, \infty) \label{eq:rv:directed:rw}
\end{align}
for any $B \in \mathscr{B}(\mathbb{R} \setminus \{0\})$ such that $\vartheta_\alpha(\partial B) = 0$. It is enough to show convergence in \eqref{eq:rv:directed:rw} for the collection of sets $\{(-\infty, -t_1) \cup (t_2, \infty) : t_1>0, t_2 > 0\}$ as these collection of intervals  is a $\pi$-system (see \cite[Lem.~6.1]{resnick:2007}).  We consider the case $B = (t, \infty)$ for $t > 0$. The set $(-\infty, t)$ with $t < 0$
can be handled similarly.
If we consider $b_n = a_n$, we get
\begin{align}
& \lim_{n \to \infty} n \prob \Big( a_n^\inv \langle \mathbf{X}, \mathbf{u} \rangle > t \Big) \nonumber \\
& = \lim_{n \to \infty} n \prob \Big( a_n^\inv \mathbf{X} \in \{ \mathbf{x} : \langle \mathbf{u}, \mathbf{x} \rangle > t \} \Big) \nonumber \\
& = t^{-\alpha} \mu \Big( \{\mathbf{x} : \langle \mathbf{u}, \mathbf{x} \rangle >1 \} \Big)
\end{align}
as $\{\mathbf{x}: \langle \mathbf{x}, \mathbf{u} \rangle > 1 \}$ is bounded away from $\mathbf{0}$ and it can be proved that $\mu$ does not put any mass at the boundary of this set. Thus, the limit exists and satisfies the scaling homogeneity property.
To complete the proof it suffices to show that $\mu \Big( \{ \mathbf{x} : \langle \mathbf{u}, \mathbf{x} \rangle > 1 \} \Big) > 0$. We show this
by using polar decomposition, invoking our assumption on the angular measure. Note that
\begin{align}
& \mu \Big( \{ \mathbf{x} : \langle \mathbf{u}, \mathbf{x} \rangle > 1 \} \Big) \nonumber \\
& = \nu_\alpha \otimes \varsigma \Big( \Big\{( r, \mathbf{y}) \in (0, \infty) \times \mathbb{S}^{d-1}: r \langle \mathbf{u}, \mathbf{y} \rangle >1  \Big\} \Big) \nonumber \\
& = \int_{\{ \mathbf{y} \in \mathbb{S}^{d-1}: \langle \mathbf{u}, \mathbf{y} \rangle > 0\}} \varsigma(d \mathbf{y}) \int_{r > (\langle \mathbf{u}, \mathbf{y} \rangle)^\inv} \nu_\alpha(dr) \nonumber \\
& = \int_{\{\mathbf{y} \in \mathbb{S}^{d-1}: \langle \mathbf{u}, \mathbf{y} \rangle > 0 \}} \Big( \langle \mathbf{u}, \mathbf{y} \rangle \Big)^{\alpha} \frac{d \varsigma}{d {\rm Leb}}(\mathbf{y}) { \rm Leb}(d \mathbf{y}).
\end{align}
It is now enough to prove that ${\rm Leb} \Big( \{ \mathbf{y} : \langle \mathbf{u}, \mathbf{y} \rangle > 0 \}  \Big) > 0$. Note that if $\mathbf{x} \in  \{\mathbf{y} \in \mathbb{S}^{d-1} : \langle \mathbf{u}, \mathbf{y} \rangle > 0 \}$, then $- \mathbf{x} \in \{ \mathbf{y} \in \mathbb{S}^{d-1}: \langle \mathbf{u}, \mathbf{y} \rangle < 0 \}$. This implies that ${\rm Leb} (\{\mathbf{y} \in \mathbb{S}^{d-1} : \langle \mathbf{u}, \mathbf{y} \rangle > 0\}) = {\rm Leb} (\{ \mathbf{y} \in \mathbb{S}^{d-1}: \langle \mathbf{u}, \mathbf{y} \rangle \neq 0\})/2$. Finally, we note that ${\rm Leb}(\{\mathbf{y} : \langle \mathbf{u}, \mathbf{y} \rangle \neq 0\}) = {\rm Leb}(\mathbb{S}^{d-1}) - {\rm Leb}(\{\mathbf{y} : \langle \mathbf{u}, \mathbf{y} \rangle =0\})$ is strictly positive, since $\{\mathbf{y} : \langle \mathbf{u}, \mathbf{y} \rangle = 0\}$ contains only $2(d-1)$ elements. Hence ${\rm Leb}(\{\mathbf{y} : \langle \mathbf{u}, \mathbf{y} \rangle =0\}) = 0$.
\end{proof}

\subsection{Proof of Theorem~\ref{thm:main:theorem:persistence}} \label{subsec:proof:persistance:one:dim}

The proof  is similar to the proof given in Section~\ref{sec:proof}. Therefore, we will provide a brief sketch of the proof below to indicate the similarity and obvious differences between these two cases.

\noindent{\bf Upper bound.} We follow the steps given in Subsection~\ref{subsec:multidimension:upper:bound}.  We follow \textbf{Step 1} with $r^* = b/a$ in the definition of $u_i$. Then the one-dimensional analogues of \eqref{eq:upper:bound:ind:increment} and \eqref{eq:upper:bound:prod:cond:prob} lead to the following inequality 
\begin{align}
& \prob \Big( S_{u_{i+1}} \in [a u_{i+1}, b u_{i+1}] \Big| S_{u_i} \in [a u_i, b u_i] \Big)  \le \prob \Big( S_{u_{i+1} - u_i} > b \eta u_i  \Big).
\end{align}
 We can again use \cite[Lemma~2.1]{hult:lindskog:mikosch:samorodnitsky:2005} with $d = 1$ to obtain the upper bound in \eqref{eq:upper:bound:final:Pn}.  Then \textbf{Step 3} produces the desired upper bound.

\noindent{\bf Lower bound.} As $\prob(X_1 \in [a, b]) > 0$, it follows that 
\begin{align*}
    \prob \Big( \bigcap_{k = 1}^N \{ k^\inv S_k \in [a, b]\} \Big) \ge \prob \Big( \bigcap_{k = 1}^N \{ X_k \in [a, b]\} \Big) = \Big[ \prob(X_1 \in [a, b]) \Big]^N > 0
\end{align*}
for any integer $N \ge 1$. Define $r = b/a$ and consider $\rho \in (0, 1 - 1/ \sqrt{r})$ which satisfies 
\begin{align}
    (1 - \rho)^2 r > 1.
\end{align}
We then define $m_i = \lfloor m_1 [(1 - \rho) r]^{i - 1} \rfloor$ for every $i \ge 2$, with $m_1$ a fixed large integer and $m_0 = 1$. We then decompose the index set $\{1, 2, \ldots, n\} = \cup_{i = 1}^{\kappa_n} D_i $ where $D_i = \{ m_i + 1, m_ i + 2, \ldots, m_{i + 1}\}$.
 We also construct a set $\Upsilon_i$ such that $ \Upsilon_i \subset  [a m_i, bm_i]$ for large enough $i$. We then enforce $S_{m_i} \in \Upsilon_i$ for all large enough $i$ which yields the lower bound to $P_n$ of the required order. By construction, we make sure that $\Upsilon_i \cap \Upsilon_{i + 1} = \emptyset$ and the distance between the sets $\Upsilon_i$ and $\Upsilon_{i + 1}$ is of the order of magnitude $m_i$. This event enforces the segment $(S_k : k \in D_i)$ \\
BZ: the notation with ":" is overloaded and it doesn't look neat. Let's discuss how to fix it. \\
to travel a distance of order $m_i$. We then write down $P_n$ in the following product form 
\begin{align}
    & \prob \Big( \bigcap_{k = 1}^{m_N - 1} \{k^\inv S_k \in [a, b]\} \cap \{ S_{m_N} \in \Upsilon_N\}  \Big) \prod_{i = N}^{\kappa_n - 1} \prob \Big( \bigcap_{j = 1}^{m_{i + 1} - m_i - 1} \nonumber \\
    & \hspace{.5cm}  \{ S_{m_i + j} \in [a (m_i + j), b(m_i + j)]\} \cap \{ S_{m_{i + 1}} \in \Upsilon_{i + 1} \} ~\Big|~ {\sf G}_i \Big) \label{eq_one_dimension_prod_formula_LB}\\
    & \mbox{ where } {\sf G}_i = \bigcap_{k = 1}^{m_N - 1} \{ k^\inv S_k \in [a, b]\} \cap \{ S_{m_N} \in \Upsilon_N\} \cap \bigcap_{j = N}^{i - 1} \Big\{ \nonumber \\
    & \hspace{.5cm} \bigcap_{j' = 1}^{m_{j + 1} - m_j - 1} \{ S_{m_j + j'} \in [a(m_j + j'), b(m_j + j')]\} \cap \{ S_{m_{j + 1}} \in \Upsilon_{j + 1}\} \Big\}.
\end{align}
By construction, the set $\Upsilon_{i + 1}$ is not accessible to the segment $(S_k : k \in D_i)$ initially. Hence, we find a positive constant $g$ such that $\Upsilon_{i+ 1}$ is accessible to $S_{\lfloor gm_i \rfloor}$ and further decompose the segment into two parts given by $(S_k : k \in D_i^{(1)})$ and $(S_k : k \in D_i^{(2)})$ where $D_i^{(1)} = \{ m_i + 1, m_i + 2, \ldots, \lfloor g m_i \rfloor\}$ and $D_i^{(2)} = \{ \lfloor g m_i \rfloor + 1, \lfloor g m_i \rfloor + 2, \ldots, m_{i + 1} \}$. In the first part of the segment, the random walk only contributes to the fluctuation (it can only travel a distance of order $O(m_i^{1/ \alpha_0 + \delta})$ where $1/\alpha_0 + \delta < 1$). The second part of the segment contains one necessary jump of order $m_i$ and the rest of the steps contribute to the fluctuation in an accumulated way.\\
To realize this strategy, we use the stationarity and the independence of the increments to write down the $i$-th term in the product formula \eqref{eq_one_dimension_prod_formula_LB} in terms of $( S_k : k \in D_i)$. The generalization of Kolmogorov's inequality (stated in \eqref{lemma:gen:kolmogorov}) is used to show that the first part $(S_k : k \in D_i^{(1)})$ can contribute to the fluctuation with high probability. The probability of the second part $(S_k : k \in D_i^{(2)})$ containing a jump of magnitude $O(m_i)$ is roughly of order $m_i^{1 - \alpha}$ leading to the right constant in Theorem~\ref{thm:main:theorem:persistence}. Thus the proof follows if we choose the constant $g$ and construct $\Upsilon_i$ in an appropriate way for large enough $i$.

We define 
\begin{align}
    & \Upsilon_i = [a(1 - \rho )r m_i + m_i^{1/\alpha_0 + \delta}, (1 - \rho/2) b m_i] \mbox{ for } i \ge 1 \nonumber \\
    & \mbox{ and } g = (1 - \rho/2)(1 - \rho) r \in (1, r(1 - \rho)).
\end{align}
To realize the strategy fully, we shall design two auxiliary sets $\widetilde{\Upsilon}_i \subseteq \cap_{j \in D_i^{(1)}} [aj, bj]$ and $\widehat{\Upsilon}_i \subseteq \cap_{j\in D_i^{(2)} \setminus \{ m_{i + 1}\}} [aj, bj]$ such that $S_k \in \widetilde{\Upsilon}_i$ for all $k \in D_i^{(1)}$ and $S_k \in \widehat{\Upsilon}_i$ for all $k \in D_i^{(2)}\setminus \{ m_{i + 1} \}$. We define 
\begin{align}
    & \widetilde{\Upsilon}_i = [a(1 - \rho) r m_i + 9 m_i^{1/\alpha_0 + \delta}, (1 - \rho/2) b m_i + m_i^{1/\alpha_0 + \delta}] \nonumber \\
    & \mbox{ and } \widehat{\Upsilon}_i = [a (1 - \rho) r m_i + 8 m_i^{1/\alpha_0 + \delta}, (1 - \rho/2) b m_{i + 1} - m_i^{1/\alpha_0 + \delta}].
\end{align}
It is easy to check that $\widetilde{\Upsilon}_i$ and $\widehat{\Upsilon}_i$ satisfy the requirements for large enough $i$ (see proof of Lemma~\ref{lemma:subsets:firstsegment:second:segment}). Therefore, we have the following lower bound on the $i$-th conditional probability in \eqref{eq_one_dimension_prod_formula_LB}:
\begin{align}
    & \prob \Big( \bigcap_{j = m_i + 1}^{\lfloor g m_i \rfloor} \{ S_j \in \widetilde{\Upsilon}_i \} ~ \Big|~ {\sf G}_i \Big) ~ \prob \Big( \bigcap_{j = \lfloor g m_i \rfloor +1}^{m_{i + 1} - 1} \{ S_j \in \widehat{\Upsilon}_i\} \nonumber \\
    & \hspace{1cm} \cap \{ S_{m_{i + 1}} \in \Upsilon_{i + 1} \} ~\Big|~ {\sf G}_i \cap \bigcap_{j = m_i + 1}^{\lfloor g m_i \rfloor} \{ S_j \in \widetilde{\Upsilon}_i \} \Big) =: {\rm T}_i^{(1)} \times {\rm T}_i^{(2)}. 
\end{align}
We shall now derive lower bounds for the terms ${\rm T}_i^{(1)}$ and ${\rm T}_i^{(2)}$ separately.

Note that the term ${\rm T}_i^{(1)}$ can be written as
\begin{align}
    \prob \Big( \bigcap_{j = m_i + 1}^{\lfloor g m_{i} \rfloor} \{ S_j \in \widetilde{\Upsilon}_i\} \cap \{ S_{m _i} \in \Upsilon_i\} \cap {\sf G}'_i \Big)/ \prob({\sf G}_i), \label{eq_one_dim_cond_prob_fraction}
\end{align}
where ${\sf G}'_i = {\sf G}_i \cup \{ S_{m_i} \notin \Upsilon_i \}$. Observe that on the event $\{ S_{m_i} \in \Upsilon_i\}$, $\{ - m_i^{1/\alpha_0 + \delta} < \min_{m_i + 1 \le j \le \lfloor g m_i \rfloor} (S_j - S_{m_i}) < \max_{m_i + 1 \le j \le \lfloor g m_i \rfloor} (S_j - S_{m_i}) < m_i^{1/\alpha_0 + \delta}\}$ implies $\{S_j \in \widetilde{\Upsilon}_i  \mbox{ for all } j \in D_i^{(1)}\}$. Therefore, we have the following lower bound for the numerator in \eqref{eq_one_dim_cond_prob_fraction}:  
\begin{align}
    & \prob \Big( {\sf G}_i \cap \Big \{ - m_i^{1/\alpha_0 + \delta} < \min_{m_i + 1 \le j \le \lfloor g m_i \rfloor} (S_j - S_{m_i}) < \nonumber \\
    & \hspace{1cm} \max_{ m_i+ 1  \le j \le \lfloor g m_i \rfloor} (S_j - S_{m_i}) < m_i^{1/\alpha_0 + \delta}  \Big\} \Big).
\end{align}
We can now use the independence of the segments $(S_k : 1 \le k \le m_i)$ and $(S_k - S_{m_i} :k \in D_i^{(1)} )$, and the distributional identity $(S_k - S_{m _i} : k \in D_i^{(1)} ) \eqd (S_j : 1 \le j \le \lfloor g m_i \rfloor - m_i)$ to obtain the following lower bound for ${\rm T}_i^{(1)}$: 
\begin{align}
    & \prob \Big( - m_i^{1/\alpha_0 + \delta} \le \min_{1 \le j \le \lfloor g m_i \rfloor - m_i } S_j \le \max_{1 \le j \le \lfloor g m_i \rfloor - m_i} S_j \le m_i^{1/\alpha_0 + \delta} \Big) \nonumber \\
    & \ge 1 - \prob \Big( \max_{1 \le j \le \lfloor g m_i \rfloor - m_i } |S_j| > m_i^{1/ \alpha_0 + \delta} \Big). \label{eq_final_lb_term_one_odim}
\end{align}
We can now use the generalized Kolmogorov's inequality when to conclude that the lower bound in \eqref{eq_final_lb_term_one_odim} is close to one if we choose $i$ large enough. 

We shall derive the exact asymptotics for the term ${\rm T}_i^{(2)}$ for large enough $i$. We want to create an envelope for the segment $(S_j : \lfloor g m_i \rfloor + 1 \le j \le m_{i + 1})$ so that the segment contains exactly one large jump (of absolute magnitude $O(m_i)$) to ensure $\cap_{j = \lfloor g m_i \rfloor + 1 }^{m_{i + 1}} \{ S_j \in \widetilde{\Upsilon}_i\} \cap \{S_{m_{i + 1}} \in \Upsilon_{i + 1}\}$. To write down the envelope explicitly, we need the following intervals
\begin{align}
    & \varpi_i = [ - m_i^{1/\alpha_0 + \delta}, (1 - \rho/2) b (m_{i + 1} - m_i) - 2 m_i^{1/\alpha_0 + \delta}], \nonumber \\
    & \widetilde{\varpi}_i = [a(1 - \rho) r (m_{i + 1} - m_i) + 10 m_{i + 1}^{1/\alpha_0 + \delta} - 9 m _i^{1/\alpha_0 + \delta}, \nonumber \\
    & \hspace{2cm}  (1- \rho/2) b (m_{i + 1} - m_i) - 3 m_i^{1/\alpha_0 + \delta}], \nonumber \\
    & \mbox{ and } \Gamma_i^* = [a(1 - \rho)r (m_{i + 1} - m_i) + 10 m_{i + 1}^{1/\alpha_0 + \delta} - 8 m_i^{1/\alpha_0 + \delta}, \nonumber \\
    & \hspace{2cm} b(1 - \rho/2) (m_{i + 1} - m_i) - 3 m_i^{1/\alpha_0 + \delta}].
\end{align}
For large enough $i$, we have the following inclusion 
\begin{align}
    & \{S_{\lfloor g m_i \rfloor} \in \widetilde{\Upsilon}_i\} \cap \Big[ \bigcap_{k = \lfloor g m_i \rfloor + 1}^{m_{i + 1} - 1} \{ S_k - S_{\lfloor g m_i \rfloor} \in \varpi_i\} \cap \{S_{m_{i + 1}} - S_{\lfloor g m_i \rfloor} \in \widetilde{\varpi}_i \} \Big] \nonumber \\
    & \hspace{1cm} \subseteq \{ S_{\lfloor g m_i \rfloor} \in \widetilde{\Upsilon}_i\} \cap \bigcap_{k = \lfloor g m_i \rfloor + 1 }^{m_{i + 1} - 1} \{ S_k \in \widehat{\Upsilon}_i \} \cap \{ S_{m_{i + 1}} \in \Upsilon_{i + 1} \}. \label{eq_one_dim_inclusion_lowerb_main}
\end{align}
We now observe that the left-hand side of the inclusion \eqref{eq_one_dim_inclusion_lowerb_main} can be decomposed into two independent events using the independent increment property of the random walk. Combining these facts, we obtain the following lower bound for the term ${\rm T}_i^{(2)}$:
\begin{align}
    \prob \Big( \cap_{j = \lfloor g m_i \rfloor + 1}^{m_{i + 1} -1} \big\{ S_j \in \varpi_i \big\} \cap \big\{ S_{m_{i + 1}} \in \widetilde{\varpi}_i \big\} \Big). \label{eq_one_dim_lb_shifted}
\end{align}
We now decompose the event inside the probability in the right hand side of \eqref{eq_one_dim_lb_shifted} into disjoint events by taking into account the location of the large jump in the interval $D_i^{(2)}$. The following event helps to write down the decomposition 
\begin{align*}
    {\sf E}_t & = \{X_t \in \Gamma_i^* \} \cap \Big\{ \big\{ \max_{1 \le k \le t - 1} S_k , \max_{t + 1 \le k \le m_{i + 1} - \lfloor g m_i \rfloor} (S_k - X_t) \big\} \le m_i^{1/\alpha_0 + \delta} \big\} \nonumber \\
    & \hspace{1cm} \cap \big\{ \min \big[ \min_{1 \le k \le t - 1} S_k, \min_{t + 1 \le k \le m_{i + 1} - \lfloor g m_i \rfloor} (S_k - X_t) \big] > - m_i^{1/\alpha_0 + \delta}  \big\} \Big\}
\end{align*}
for every $t \in D_i^{(2)}$. It is easy to check that $\cup_{1 \le t \le m_{i + 1} - \lfloor g m_i \rfloor} {\sf E}_t$ implies the event inside the probability in  \eqref{eq_one_dim_lb_shifted}. We can now use exchangeability of the random variables $(X_t : 1 \le t \le m_{i + 1} - \lfloor g m_i \rfloor)$ to see that $\prob({\sf E}_t) = \prob({\sf E}_1)$ for every $t \ge 1$ and obtain the following lower bound for ${\rm T}_i^{(2)}$:
\begin{align}
    & (m_{i + 1} - \lfloor g m_i \rfloor) \prob({\sf E}_1) \nonumber \\
    & = (m_{i + 1} - \lfloor g m_i \rfloor) \prob \big( X_1 \in \Gamma_i^* \big) \prob \Big( \Big\{ \max_{1 \le k \le m_{i + 1} - \lfloor g m_i \rfloor - 1} S_k \le m_i^{1/\alpha_0 + \delta} \Big\} \nonumber \\
    & \hspace{1cm} \cap \Big\{ \min_{1 \le k \le m_{i + 1} - \lfloor g m_i \rfloor - 1 } \ge - m_i^{1/\alpha_0 + \delta} \Big\} \Big). \label{eq_one_dim_lb_location_large_jump}
\end{align}
For large enough $i$, the last probability in \eqref{eq_one_dim_lb_location_large_jump} is very close to $1$ as we have seen earlier in the analysis of term ${\rm T}_i^{(1)}$ and so, we can ignore that for the further analysis.  We can use now regular variation to conclude that 
\begin{align}
    & (m_{i + 1} - \lfloor g m_i \rfloor) \prob \big( X_1 \in \Gamma_i^* \big) \nonumber \\
    & \sim   \big[ (1 - \rho) r - g \big] m_i \prob \Big( m_i^\inv X_1 \in [a(1 - \rho)^2 r^2, b(1 - \rho/2)(1 - \rho) r] \Big) \nonumber \\
    & \sim [(1 - \rho) r - g] \big [ \alpha \int_{a(1 - \rho)^2 r ^2}^{b(1 - \rho/2)(1 - \rho)} x^{- \alpha - 1} \dtv x  \big]m_i^{1 - \alpha} \nonumber \\
    & \sim {\rm const.} ~ \exp \big\{ - i \big[ (\alpha - 1) \log \big( (1 - \rho) r \big)  \big]  \big\},
\end{align}
as $i \to \infty$. The lower bound now follows from simple algebra (see \eqref{eq_upper_bound_prodct_formula} in Step~3 in the proof of \eqref{eq:upper:bound:aim}), and by letting $\rho \to 0$.

\section*{Acknowledgement}
The authors are thankful to 
Guido Janssen for an analytic computation which led to a correct guess of the proper normalization of $P_n$.
The authors are thankful to  a referee for valuable suggestions which improved the quality of the exposition.


\bibliographystyle{plain}

\end{document}